\magnification 1200
\let \ol=\overline
\let \ti=\widetilde

\font \srm=cmr10 at 7.5pt
\font \ssrm=cmr10 at 5.625pt
\font\ens=msbm10
\font\sens=msbm10 at 7.5pt
\font\main=cmsy10 at 10pt
\font\smain=cmsy10 at 7.5pt
\font\ssmain=cmsy10 at 5.625pt
\font \para=cmr10 at 15.4 pt

\font \fin=lasy8 at 15.4 pt
\def \o{\mathop{\hbox{\main O}}\nolimits}
\def \so{\mathop{\hbox{\smain O}}\nolimits}
\def \sso{\mathop{\hbox{\ssmain O}}\nolimits}
\def \mes{\mathop{\hbox{mes}}\nolimits}
\def \H{\mathop{\hbox{\rm PW}}\nolimits}

\def \P{\mathop{\hbox{\main P}}\nolimits}

\def \sP{\mathop{\hbox{\smain P}}\nolimits}

\def \End{\mathop{\hbox{\rm End}}\nolimits}
\def \Hom{\mathop{\hbox{\rm Hom}}\nolimits}
\def \sp{\mathop{\hbox{\rm sp}}\nolimits}
\def \hors{\mathop{\hbox{\rm hors}}\nolimits}
\def \supp{\mathop{\hbox{\rm supp}}\nolimits}
\centerline {\para Une formule de Plancherel pour l'alg\`ebre de Hecke}

\null
\centerline {\para d'un groupe r\'eductif p-adique}

\null
\null
\centerline {Volker Heiermann}

\centerline {Institut f\"ur Mathematik, Humboldt-Universit\"at zu Berlin}
\centerline {D-10099 Berlin}
\centerline {$e$-mail: heierman@mathematik.hu-berlin.de}

\null
\null
\null
{\bf R\'esum\'e:} Nous montrons un th\'eor\`eme de Paley-Wiener matriciel
pour l'alg\`ebre de Hecke d'un groupe r\'eductif $p$-adique. La preuve est bas\'ee sur une analogue de la formule de Plancherel.

\null
\null
\null
Fixons un corps local non archim\'edien $F$. Soit $G$ l'ensemble des points
$F$-rationnels d'un groupe r\'eductif connexe $\underline {G}$ d\'efini sur $F$. Fixons un sous-groupe ouvert compact $K$ maximal sp\'ecial de $G$. On munit tout sous-groupe alg\'ebrique ferm\'e $H$ de $G$ de la mesure de Haar invariante \`a gauche pour laquelle mes$(H\cap K)=1$. Lorsque $M$ est un sous-groupe de L\'evi de $G$ (ou plus pr\'ecis\'ement l'ensemble des points $F$-rationnels d'un facteur de L\'evi d'un sous-groupe parabolique de $\underline {G}$ d\'efini sur $F$), notons $X(M)$ le groupe des caract\`eres non ramifi\'es de $M$ (d\'efini en 1.2). C'est une vari\'et\'e alg\'ebrique complexe isomorphe \`a $(\hbox{\ens C}^{\times })^r$, o\`u $r$ d\'esigne la dimension du tore d\'eploy\'e maximal dans le centre de $M$. Pour une repr\'esentation cuspidale irr\'eductible $(\sigma ,E)$ de $M$, on notera $\o _{\sigma }=\{\sigma \otimes \chi \vert \chi \in X(M)\}$ son \it orbite inertielle. \rm L'application
$X(M)\rightarrow \o _{\sigma }$, $\chi \mapsto \sigma \otimes \chi$, d\'efinit de fa\c con naturelle une structure de vari\'et\'e alg\'ebrique complexe sur $\o _{\sigma }$. Une fonction complexe $\varphi $ sur $\o _{\sigma }$ sera dite polynomiale (resp. rationnelle), si la fonction $\chi\mapsto \varphi (\sigma \otimes \chi )$ est polynomiale (resp. rationnelle) sur $X(M)$.

Lorsque $P$ est un sous-groupe parabolique de L\'evi $M$, on d\'esigne par $i_P^G$ le foncteur de l'induction parabolique unitaire. Si $M$ et $K$ sont en bonne position relative, on d\'efinit l'espace $i_{P\cap K}^KE$ des applications $f:K\rightarrow E$ invariantes \`a droite par un sous-groupe ouvert et v\'erifiant $f(muk)=\sigma (m)f(k)$ pour tout $m\in M\cap K, u\in U\cap K$ et $k \in K$. La restriction $i_P^G E \rightarrow i_{P\cap K}^K E$ est un isomorphisme, et l'espace $i_{P\cap K}^K E$ ne change pas si on remplace $\sigma $ par un \'el\'ement de son orbite inertielle $\o _{\sigma }$. Toutes les
repr\'esentations $i_P^G\sigma '$, $\sigma '\in \o _{\sigma }$, se r\'ealisent donc dans le m\^eme espace $i_{P\cap K}^K E$. Ceci permet d'introduire la notion naturelle d'une application polynomiale sur $\o _{\sigma }$ \`a valeurs dans $\Hom (i_{P\cap K}^KE, i_{P'\cap K}^KE)$ ou dans $i_{P\cap K}^KE\otimes i_{P\cap K}^KE^{\vee }$ (cf. [W], IV.1 et VI.1), $P'$ d\'esignant un deuxi\`eme sous-groupe parabolique de L\'evi $M$.

L'op\'erateur d'entrelacement $J_{P'\vert P}(\sigma ')$ est d\'efini pour $\sigma '$ dans un certain ouvert Zariski dense de $\o _{\sigma }$. C'est une application lin\'eaire de $i_{P\cap K}^KE$\ dans $i_{P'\cap K}^K E$ qui v\'erifie $J_{P'\vert P}(\sigma ')$ $(i_P^G \sigma ')(g)=(i_{P'}^G\sigma ')(g) J_{P'\vert P}(\sigma ')$ pour tout $g\in G$, $\sigma '\in \o _{\sigma }$. On a
$\displaystyle {<(J_{P'\vert P}(\sigma ')v)(g),e^{\vee }>}$ $\displaystyle {=\int _{U\cap U'\backslash U'} <v(u'g),e^{\vee }> du'}$ pour $v\in i_{P\cap K}^K E$ et $e^{\vee }\in E^{\vee }$, si l'int\'egrale \`a droite est convergente. L'application $\o _{\sigma }\rightarrow \hbox{Hom}(i_{P\cap K}^K E,i_{P'\cap K}^K E), \sigma '\mapsto J_{P'\vert P}(\sigma ')$, est rationnelle (i.e. il existe une fonction polynomiale $p$ sur $\o _{\sigma }$ telle que l'application $\sigma '\mapsto p(\sigma ')J_{P'\vert P}(\sigma ')$ soit polynomiale sur $\o _{\sigma }$). Pour la preuve de ces r\'esultats et d'autres propri\'et\'es des op\'erateurs d'entrelacement, nous renvoyons le lecteur \`a [W]. Remarquons que la plupart des r\'esultats qui y sont expos\'es sont dus \`a Harish-Chandra.

Fixons un tore d\'eploy\'e maximal $A_0$ de $G$ par rapport auquel $K$ est en bonne position et notons $W^G:=W(G,A_0)$ le groupe de Weyl d\'efini relatif \`a ce tore. Si $M$ est semi-standard (i.e. $M\supseteq A_0$) et que $w\in W^G$, on dispose d'un isomorphisme $\lambda (w): i_P^G E\rightarrow i_{wP}^G wE, v\mapsto v_w, v_w(g):=v(w^{-1}g)$, entre les repr\'esentations $i_P^G\sigma $ et $i_{wP}^G (w\sigma )$.

Notons $C_c^{\infty }(G)$ l'espace des fonctions complexes lisses \`a support compact sur $G$. Il est bien connu que l'on peut associer \`a tout \'el\'ement $f$ de $C_c^{\infty }(G)$ un endomorphisme $(i_P^G\sigma )(f)$ de l'espace vectoriel $i_P^G E$ que nous noterons $\widehat {f}^G(P,\sigma )$. Notre but est le r\'esultat suivant:

\null
{(\bf 0.1)} {\bf Th\'eor\`eme:} \it
\'Etant donn\'e pour chaque (classe d'\'equivalence d'une)
repr\'esen-tation cuspidale irr\'eductible $(\sigma, E_{\sigma })$ d'un sous-groupe de L\'evi
semi-standard $M$ de $G$, et tout sous-groupe parabolique $P$ de L\'evi $M$,  un endomorphisme
$\varphi _{P,\sigma }$ de l'espace vectoriel $i_{P\cap K}^K E_{\sigma }$ tels que la famille
$\{\varphi _{P,\sigma }\}_{(P,\sigma )}$ v\'erifie les propri\'et\'es suivantes:

1) Pour tout $(P,\sigma )$, l'application $\varphi _{P,\so }:\o \rightarrow
\hbox{End}(i_{P\cap K}^K E_{\sigma })$, $\sigma '\mapsto \varphi _{P,\sigma '}$, est polynomiale
sur l'orbite inertielle $\o$ de $\sigma $;

2) Il existe un sous-groupe ouvert compact de $G$ par lequel toute composante $\varphi _{P,\sigma }$ 
est invariante \`a gauche et \`a droite;

3) Pour tout $(P,\sigma )$ et tout $w\in W^G$, on a $\lambda (w)\circ \varphi _{P,\sigma
}=\varphi _{wPw^{-1},w\sigma }\circ \lambda (w)$;

4) Pour tout $(P,\sigma )$ et tout $(P',\sigma )$, on a l'identit\'e d'applications rationnelles $J_{P'\vert P}(\sigma )\circ \varphi
_{P,\sigma }=\varphi _{P',\sigma }\circ J_{P'\vert P}(\sigma )$;

{\noindent alors il existe une fonction $f$ dans $C_c^{\infty }(G)$, telle que $\varphi
_{P,\sigma }=\widehat {f}^G(P,\sigma )$ pour tout $(P,\sigma )$.}\rm

\null
R\'eciproquement, il est bien connu que, pour $f$ dans $C_c^{\infty }(G)$, la famille $\{\widehat {f}^G(P,\sigma )\}_{(P,\sigma )}$ v\'erifie les propri\'et\'es 1) - 4) du th\'eor\`eme 0.1.

La propri\'et\'e 2) \'equivaut \`a dire que l'image de $\varphi _{P,\so }$ est inclus dans un sous-espace de dimension finie de l'image de l'application canonique $i_P^GE\otimes i_P^GE^{\vee }\hookrightarrow \End (i_P^GE)$, et qu'il n'existe qu'un nombre fini de $(P,\o )$ avec $\varphi _{P,\so }\ne 0$ (cf. [W] th\'eor\`eme VIII.1.2). 

Notre d\'emonstration de ce th\'eor\`eme est bas\'ee sur une analogue de la formule de Plancherel de Harish-Chandra. Elle a donc l'avantage d'expliciter la fonction $f$ du th\'eo-r\`eme. La preuve utilise le r\'esultat suivant qui sera prouv\'e dans la partie B:

\null
{\bf (0.2)} {\bf Proposition:} \it
Soit $\{\varphi _{P,\so }\}_{(P,\so )}$ comme dans le th\'eor\`eme. Pour tout $(P,\o )$, il existe une application polynomiale $\xi _{P,\so }:\o \rightarrow i_{\overline {P}\cap K}^KE_{\so }\otimes i_{P\cap K}^K E_{\so }^{\vee }$ (o\`u $E_{\so }:=E_{\sigma }$ pour un $\sigma
\in \o$) \`a image dans un espace de dimension finie, telle que
$$\varphi _{P,\so }(\sigma )=\sum _{w\in W; w\so =\so } (J_{P\vert \overline {wP}}(\sigma )\circ\lambda (w))\otimes (J_{P\vert wP}(\sigma ^{\vee })\circ \lambda (w))\ \xi _{P,\so }(w^{-1}\sigma ),$$
pour tout $\sigma \in \o $.\rm

\null
{\noindent (Ici on a identifi\'e $\varphi _{P,\so }(\sigma )\in \hbox{End}(i_{P\cap K}^K E)$
\`a un \'el\'ement de $i_{P\cap K}^K E\otimes i_{P\cap K}^K E^{\vee }$.)}
\null

Plus pr\'ecis\'ement, choisissons $\xi _{P,\so }$ et posons $\zeta _{P,\so }(\sigma )=(J_{\overline {P}\vert P}(\sigma )^{-1}\otimes 1) \xi _{P,\so }(\sigma )$. Pour $\sigma\in\o $, notons $E_{P,\sigma }^G$ l'application lin\'eaire qui associe \`a un \'el\'ement $v\otimes v^{\vee }$
de $i_{P\cap K}^K E_{\so }\otimes i_{P\cap K}^K E_{\so }^{\vee }$ la fonction $g\mapsto $
$<(i_P^G\sigma )(g)v,v^{\vee }>$, $g\in G$. Fixons $\sigma \in \o$. Avec $\varphi =\varphi _{P,\so }$, posons
$$f_{\varphi }(g)=\int _{\hbox {\srm Re}(\chi )=\mu >>_P 0} E_{P,\sigma \otimes \chi }^G(\zeta _{P,\so }(\sigma \otimes \chi ))(g^{-1})\ dIm(\chi ),$$
pour $g\in G$. (La partie r\'eelle d'un caract\`ere non ramifi\'e \'etant d\'efinie dans 1.2, la notation $\mu >>_P 0$ est justifi\'ee par le fait que l'on peut trouver $\mu $ dans la chambre de Weyl de $P$ tel que les p\^oles $\chi $ de la fonction dans l'int\'egrale v\'erifient $<\alpha ^{\vee },\hbox {\rm Re}(\chi )>\ <\ <\alpha ^{\vee },\mu >$ pour toute racine $\alpha $ positive pour $P$. Gr\^ace au th\'eor\`eme des r\'esidus, la valeur de l'int\'egrale ne d\'epend alors pas du choix de $\mu $ v\'erifiant cette condition.)

On montrera (cf. A.3.1-2):

1) La fonction $f_{\varphi _{P,\sso }}$ ne d\'epend que de $\varphi _{P,\so }$;

2) La fonction $f_{\varphi _{P,\sso }}$ appartient \`a $C_c^{\infty }(G)$;

Notons $[\o]=\{w\sigma \vert  \sigma \in \o ,w \in W^G\}$ la classe inertielle de $\o $ et $M$ le sous-groupe de L\'evi semi-standard sous-jacent \`a $\o $. Posons $f_{\varphi _{[\sso ]}}=c([\o ])\sum _{(P',\so ')} f_{\varphi _{P',\sso '}}$, o\`u $c([\o ])$ est une constante pr\'ecis\'ee dans 3.2, la somme portant sur les couples $(P',\o ')$ form\'es d'une orbite inertielle $\o '=w\o$, $w\in W^G$, et d'un sous-groupe parabolique de L\'evi $wMw^{-1}$. Alors on a par ailleurs

3) $\widehat {f_{\varphi _{[\sso ]}}}(P,\sigma )=0$ si $\sigma \not\in [\o ]$;

4) $\widehat {f_{\varphi _{[\sso ]}}}(P,\sigma )=\varphi _{P,\sigma }$ si $\sigma \in [\o ]$;

5) $\int _G f_{\varphi _{[\sso ]}}(g)\ \overline {f_{\varphi _{[\sso ']}}(g)} dg=0$ si $[\o
]\ne [\o' ]$;

6) La fonction $f_{\varphi }=\sum _{[\sso ]} f_{\varphi _{[\sso ]}}$ v\'erifie $\varphi _{P,\sigma }=\widehat {f}(P,\sigma )$ pour tout $(P,\sigma )$.

\null
Cet article est divis\'e en deux parties. Dans la partie A, nous prouvons tous les r\'esultats annonc\'es dans l'introduction \`a l'exception de la proposition 0.2. Sa preuve est le contenu de la partie B. Les deux parties peuvent \^etre lues ind\'ependamment, seules certaines d\'efinitions et notations introduites dans la section A.1 seront utilis\'ees sans rappel dans la partie B.

\null
Remarquons que J.N. Bernstein a annonc\'e une preuve du th\'eor\`eme 0.1 par une m\'ethode diff\'erente de la n\^otre.

\null
L'essentiel de ce travail a \'et\'e r\'ealis\'e, lorsque l'auteur s\'ejournait \`a l'Universit\'e Paris 7 au sein de l'\'equipe "Th\'eorie des Groupes". Ce s\'ejour a \'et\'e financ\'e par une bourse Feodor Lynen de la fondation Alexander von Humboldt en correspondance avec M.-F. Vign\'eras. Cette bourse comprenait une participation financi\`ere de l'Universit\'e Paris 7 venant du r\'eseau "G\'eom\'etrie arithm\'etique alg\'ebrique" soutenu par le programme "Formation et Mobilit\'e des Chercheurs" de l'Union Europ\'eenne. Mon tuteur aupr\`es de la fondation Alexander von Humboldt \'etait E.-W. Zink.

\null
Mes remerciements vont par ailleurs tout particuli\`erement \`a J.-L. Waldspurger \`a qui je dois l'id\'ee pour ce travail et qui m'a bien accompagn\'e durant sa r\'ealisation.

\null
\null
\null
{\bf A. Une analogue de la formule de Plancherel}

\null
\null
\null
{\bf 1.} On garde les notations et d\'efinitions de l'introduction. On notera $q$ le cardinal du corps r\'esiduel de $F$, $v_F$ la valuation discr\`ete de $F$ normalis\'ee telle que $v_F(F^{\times })=\hbox{\ens Z}$ et $\vert \vert _F$ la valeur absolue, donn\'ee par $\vert x \vert _F=q^{-v_F(x)}$ pour $x$ dans $F^{\times }$.

Les d\'efinitions et notations qui seront introduites dans la suite pour le groupe r\'eductif $G$ et munies du symb\^ole $G$ seront ensuite utilis\'ees pour tout groupe r\'eductif $M$, en rempla\c cant $G$ par $M$, sans que cela soit dit explicitement.

\null
{\bf 1.1} Notons $A_G$ le tore d\'eploy\'e maximal dans le centre de $G$ et $G^{\hbox {\srm der}}$ le groupe d\'eriv\'e de $G$. Posons $X^*(G)=\Hom _F(G,\hbox {\ens G}_m)$ et $X_*(S)=\Hom _F(\hbox {\ens G}_m,S)$ lorsque $S$ est un tore.

\null
{\bf 1.1.1} Posons $a_0^*=X^*(A_0)\otimes _{\hbox {\sens Z}}\hbox {\ens R}$, $a_G^*=X^*(A_G)\otimes _{\hbox{\sens Z}}\hbox {\ens R}$ et $a_0^{G*}=X^*(A_0\cap G^{\hbox {\srm der}})\otimes _{\hbox{\sens Z}}\hbox {\ens R}$. On a une d\'ecomposition $a_0^*=a_G^*\oplus a_0^{G*}$. Lorsque $S$ est un tore d\'eploy\'e dans $A_0$, on notera $\Sigma (S)$ l'ensemble des racines pour l'action adjointe de $S$ dans l'alg\`ebre de Lie de $G$. Soit $(P,M)$ un couple parabolique semi-standard, i.e. $P$ est un sous-groupe parabolique de $G$, $M$ un facteur de L\'evi de $P$, et on a $M\supseteq A_0$. L'ensemble $\Sigma (A_M)$ est l'ensemble des projections non nulles dans $a_M^*$ d'\'el\'ements de $\Sigma (A_0)$ suivant la d\'ecomposition $a_0^*=a_M^*\oplus a_0^{M*}$. On notera $\Sigma (P)$ l'ensemble des racines pour l'action adjointe de $A_M$ dans l'alg\`ebre de Lie du radical unipotent de $P$. 

\null
{\bf 1.1.2} On fixera pour la suite un couple parabolique semi-standard $(P_0,M_0)$ avec $P_0$ minimal. On a alors $M_0=Z_G(A_0)$ et $A_{M_0}=A_0$. L'ensemble $\Sigma (A_0)$ est un syst\`eme de racines dans $a_0^{G*}$. Remarquons que ce syst\`eme de racines peut ne pas \^etre r\'eduit. Les \'el\'ements de $\Sigma (P_0)$ s'identifient aux racines positives dans $\Sigma (A_0)$ pour un certain ordre sur $a_0^{G*}$. La base de $\Sigma (A_0)$ correspondant \`a cet ordre sera not\'e $\Delta $. Un couple parabolique $(P,M)$ sera dit standard, s'il est semi-standard et $P\supseteq P_0$. On a une bijection $\Omega \mapsto (P_{\Omega },M_{\Omega })$ entre les sous-ensembles de $\Delta $ et les couples paraboliques standards, les racines dans $\Sigma (A_0)$ de restriction triviale \`a $A_{M_{\Omega }}$ \'etant les combinaisons lin\'eaires de $\Omega $.

\null
{\bf 1.1.3} Posons $a_0=X_*(A_0)\otimes _{\hbox{\sens Z}}\hbox {\ens R}$, $a_G=X_*(A_G)\otimes _{\hbox {\sens Z}}\hbox {\ens R}$ et $a_0^{G}=X_*(A_0\cap G^{\hbox {\srm der}})\otimes_{\hbox{\sens Z}}\hbox {\ens R}$. Les espaces $a_0$ et $a_0^*$ sont duaux, l'orthogonal de $a_G$ dans $a_0^*$ est $a_0^{G*}$ et celui de $a_0^G$ est $a_G^*$.

\null
{\bf 1.1.4} Si $M$ est un sous-groupe de L\'evi semi-standard, il existe une notion de coracine $\alpha ^{\vee }$ associ\'ee \`a une racine $\alpha\in\Sigma (A_M)$. C'est un \'el\'ement de $a_M$. On en d\'eduit, pour tout sous-groupe parabolique $P$ de L\'evi $M$, une notion de chambre de Weyl dans $a_M^*$ qui est l'ensemble des \'el\'ements positifs pour $P$. 

\null
{\bf 1.1.5} On d\'efinit une application $H_0:M_0\rightarrow \Hom (X^*(M_0), \hbox{\ens R})\simeq a_0$ par $<\chi , H_0(m)>\ =v_F(\chi (m))$. Soit $(P,M)$ un couple parabolique semi-standard. Un \'el\'ement $a\in A_M$ sera dit positif pour $P$, si $<\alpha ,H_0(m)>\geq 0$ pour tout $\alpha \in \Sigma (P)$. On dira qu'il est strictement positif, si l'on a l'in\'egalit\'e stricte pour tout $\alpha \in \Sigma (P)$.
 
\null
\null
{\bf 1.2} La restriction $X^*(G)\rightarrow X^*(A_G)$ induit un isomorphisme $X^*(G)\otimes _{\hbox {\sens Z}}\hbox {\ens R}\rightarrow a_G^*$. Le groupe $X (G)$ des caract\`eres non ramifi\'es de $G$ est par d\'efinition l'image de l'homomor-phisme $a_{G,\hbox {\sens C}}^*=a_G^*\otimes _{\hbox {\sens R}}
{\hbox {\ens C}}\rightarrow \Hom(G, \hbox{\ens C}^{\times })$ qui associe \`a $\lambda=\alpha \otimes s$ le caract\`ere $\chi _{\lambda }$ tel que $\chi_{\lambda }(g)=\vert\alpha (g)\vert _F^s$. Son noyau est de la forme ${2\pi i\over \log q}R_G$, o\`u $R_G$ d\'esigne un r\'eseau de rang maximal dans 
$X^*(G)\otimes_{\hbox {\sens Z}}\hbox{\ens Q}$. L'homomorphisme munit $X (G)$ d'une structure de vari\'et\'e alg\'ebrique complexe isomorphe \`a $(\hbox{\ens C}^{\times })^r$ avec $r=$rang de $A_G$. Sa restriction \`a $a_G^*$ induit un isomorphisme avec le sous-groupe des caract\`eres r\'eels \`a valeurs $>0$. La partie r\'eelle d'un caract\`ere non ramifi\'e $\chi $, not\'e Re$(\chi )$, est l'unique \'el\'ement $\lambda $ de $a_G^*$ qui v\'erifie $\chi _{\lambda }=\vert \chi \vert $. On notera $X _{\hbox {\srm im}} (G)$ le sous-groupe de $X (G)$ form\'e des $\chi $ avec Re$(\chi )=0$. 

\null
{\bf 1.2.1} On munit $X _{\hbox {\srm im}} (A_G)$ de la mesure de Haar de masse totale $1$, et $X _{\hbox {\srm im}} (G)$ de la mesure de Haar pour laquelle la restriction $X _{\hbox {\srm im}} (G)\rightarrow X _{\hbox {\srm im}} (A_G)$ pr\'eserve localement les mesures. 
 
\null
{\bf 1.2.2} Soit $(P,M)$ un couple parabolique semi-standard. Soit $r$ une fonction rationnelle sur $X (M)$. Supposons qu'il existe un nombre fini de hyperplans de la forme $<\lambda ,\alpha ^{\vee }>=c$ dans $a_M^*$, $\alpha \in \Sigma (P)$, tels que tout p\^ole $\chi $ de $r$ soit de la forme $\chi =\chi _{\lambda }$ avec $\lambda $ sur un de ces hyperplans. Il r\'esulte du th\'eor\`eme des r\'esidus que l'int\'egrale $\int _{X_{\hbox {\ssrm im}} (M)} r(\chi _{_0}\chi ) d\chi $ reste constante, si $\chi _{_0}$ varie dans l'ouvert de $X (M)$ d\'efini par les in\'egalit\'es $<\hbox{\rm Re}(\chi ),\alpha ^{\vee }> < \ <\hbox {Re}(\chi _{_0}), \alpha ^{\vee }>$, $\chi $ parcourant les p\^oles de $r$, $\alpha \in \Sigma (P)$.

On \'ecrira plus simplement $\int _{\hbox {\srm Re}(\chi )=\mu >>_P0} r(\chi )d(Im(\chi
))$ pour la valeur de cette int\'e-grale.

L'expression $\int _{\hbox {\srm Re}(\chi )=\mu <<_P0} r(\chi )d(Im(\chi
))$ aura la signification \'evidente.

\null
{\bf 1.2.3 Proposition:} \it
Soient $D$ un ouvert de $a_G^*$ et $\psi $ une fonction holomorphe dans l'ouvert de $X (G)$ form\'e des points $\chi $ avec \hbox {\rm Re}$(\chi )\in D$. Fixons $\mu \in D$.

Alors, pour tout $\chi _{_0}\in X (G)$, $\hbox {\rm Re}(\chi _{_0})=\mu $, on 	a
$$\sum _{a\in A_G\cap K\backslash A_G} \chi _{_0}^{-1}(a) \int _{\hbox {\srm Re}(\chi )=\mu }
\psi (\chi )\chi (a) dIm(\chi)\ =\sum _{\chi } \psi (\chi \chi _{_0}),$$
la somme portant sur les \'el\'ement de $X _{\hbox {\srm im}} (G)$ de restriction triviale \`a $A_G$.

\null
Preuve: \rm Ceci r\'esulte de la th\'eorie de Fourier sur un tore. \hfill {\fin 2}

\null
\null
{\bf 1.3} Fixons un sous-groupe de L\'evi semi-standard $M$ de $G$. Notons {\main P}$(M)$ l'ensemble des sous-groupes paraboliques $P$ de $G$ de la forme $P=MU$. Fixons une repr\'esentation irr\'eductible cuspidale $(\sigma ,E)$ de $M$. Soient $P, P' \in \hbox{\main P}(M)$. Les points, o\`u l'application rationnelle $X (M)\rightarrow \Hom (i_{P\cap K}^KE, i_{P'\cap K}^KE)$, $\chi \mapsto J_{P'\vert P}(\sigma \otimes \chi )$, a un p\^ole ou bien o\`u $J_{P'\vert P}(\sigma \otimes \chi )$ n'est pas inversible, sont de la forme $\chi =\chi _{\lambda }$ avec $\lambda $ sur un nombre fini de hyperplans  de $a_M^*$ de la forme $<\alpha ^{\vee },\lambda >=c$, $\alpha \in \Sigma (P')\cap \Sigma (\overline {P})$. Soit $P''\in \P(M)$. Il existe une fonction rationnelle $j_{P\vert P'\vert P''}$ sur l'orbite inertielle $\o $ de $\sigma $, telle que $J_{P\vert P'}(\sigma\otimes \chi )J_{P'\vert P''}(\sigma \otimes \chi )=j_{P\vert P'\vert P''}(\sigma \otimes \chi )J_{P\vert P''}(\sigma \otimes \chi )$ pour tout $\chi $. Si $P''=P$, on \'ecrira plus simplement $j_{P\vert P'}=j_{P\vert P'\vert P''}$. L'\'egalit\'e $j_{P\vert P'\vert P''}(\sigma \otimes \chi )=1$ vaut si $d(P\vert P'')=d(P\vert P')+d(P'\vert P'')$,  $d(P\vert P'')$ d\'esignant le nombre d'hyperplans s\'eparant les chambres de Weyl de $a_M^*$ correspondant \`a $P$ et $P''$ respectivement. 

\null
{\bf 1.3.1} Le th\'eor\`eme suivant est un r\'esultat cl\'e pour la suite. (A notre connaissance il est paru pour la premi\`ere fois (en tous cas le cas temp\'er\'e) dans les papiers de Casselman et dans ceux de Harish-Chandra, sans que nous nous sentions comp\'etent de l'attribuer \`a l'un ou l'autre. 
Il nous a sembl\'e que ce que nous faisons rel\`eve davantage  de Casselman, Harish-Chandra adoptant un point de vue tr\`es analytique.)

Soit $(P',M')$ un couple parabolique semi-standard. Posons $W(M,M')=W^{M'}\backslash \{w\in W^G\vert \ wMw^{-1} \subseteq M'\}$ et identifions ses \'el\'ements \`a certains \'el\'ements de $W^G$. Les formules qui suivent seront essentiellement ind\'ependantes du choix d'une telle identification. Pour $w\in W(M,M')$, d\'efinissons $P_w', \ti {P}_w'\in \hbox {\main P}(M)$ par $P_w'=(w^{-1}M'w\cap P)w^{-1}U'w$ et $\ti{P}_w'=(w^{-1}M'w\cap P)w^{-1}\overline {U'}w$.

On associe par ailleurs \`a $M'$ comme dans [W] une constante $\gamma (G\vert M')$.

Pour $v\in i_{P\cap K}^KE$, $v^{\vee }\in i_{P\cap K}^KE^{\vee }$, $\chi \in X (M)$ et $a\in A_{M'}$, posons $c_{P'\vert P}(\sigma \otimes \chi ,w)(v\otimes v^{\vee })(a)= <i_{wP\cap M'}^{M'}(w(\sigma \otimes \chi ))(a)(\lambda (w)J_{P_w'\vert P}(\sigma\otimes\chi ) v)_{\vert M'},$
$(\lambda (w)J_{\ti {P}_w'\vert P}(\sigma ^{\vee }\otimes \chi ^{-1})v^{\vee })_{\vert M'}>_{M'}$.

\null
{\bf Th\'eor\`eme:} \it
Soient $v\in i_{P\cap K}^KE$ et $v^{\vee }\in i_{P\cap K}^KE^{\vee }$. Il existe $t>0$, de sorte que pour tout $\chi \in X (M)$ et tout $a\in A_{M'}$ tel que $<\alpha ,H_0(a)>\ >t$ pour tout $\alpha \in \Sigma (P')$, on ait 
$$\eqalign {
<i_P^G(\sigma \otimes \chi )(a)v, v^{\vee }>
=\gamma (G\vert M')^{-1}\delta _{P'}^{1/2}(a)\sum _{w\in W(M,M')} 
c_{P'\vert P}(\sigma \otimes \chi ,w)(v\otimes v^{\vee })(a).
}$$

\null
Preuve: \rm
La preuve du cas temp\'er\'e est faite dans [W] dans la d\'emonstration de la proposition V.1.1. Le principe de cette preuve se g\'en\'eralise \`a notre situation: il suffit de remplacer le calcul du produit bilin\'eaire de Casselman $<,>_{P'}$ relatif au module de Jacquet \it faible \rm par celui relatif au module de Jacquet. Le calcul suit les m\^emes lignes. \hfill {\fin2 }

\null
{\bf 1.3.2} Le r\'esultat suivant sera utile lors des applications du th\'eor\`eme 1.3.1.

\null
{\bf Lemme:} \it Soit $w\in W(M,M')$. On a
$$J_{P_w'\vert P}(\sigma ')J_{\overline {P}\vert P}(\sigma ')^{-1}
=j_{P_w'\vert \overline {P}}(\sigma ')^{-1} J_{P_w'\vert \overline {P}}(\sigma ')$$
en tout point $\sigma '$ de l'orbite inertielle de $\sigma $ en lequel ces op\'erateurs
sont d\'efinis.

\null
Preuve: \rm
Par la formule du produit (cf. [W] p. 55), on a
$$J_{P_w'\vert \overline {P}}(\sigma ')J_{\overline {P}\vert P}(\sigma ')=
j_{P_w'\vert \overline {P}}(\sigma ')J_{P_w'\vert P}(\sigma ').$$
\hfill {\fin 2}

\null
\null
{\bf 1.4.} Soient $(\pi ,V)$ et $(\pi ',V')$ deux repr\'esentations irr\'eductibles cuspidales de $G$. Supposons que les restrictions \`a $A_G$ de leurs caract\`eres centraux co\"\i ncident. Pour $v\in V$, $v^{\vee }\in V^{\vee }$, $v'\in V'$ et $v'^{\vee }\in V'^{\vee }$, posons
$$I(v,v^{\vee },v',v'^{\vee })=\int _{A_G\backslash G}
<\pi (g) v,v^{\vee }>\ <v',\pi '^{\vee }(g)v'^{\vee
}>dg.$$

\null
{\bf Th\'eor\`eme:} \it
i) Si $\pi \not\simeq \pi '$, alors $I(v,v^{\vee },v',v'^{\vee })=0$ pour tous $v$, $v^{\vee }$,
$v'$ et $v'^{\vee }$.

ii) Si $(\pi ,V)=(\pi ',V')$, il existe un r\'eel $d(\pi )>0$, appel\'e le degr\'e formel de $\pi$,
tel que $I(v,v^{\vee },v',v'^{\vee })=d(\pi )^{-1} <v,v'^{\vee }> <v',v^{\vee }>$ pour tous $v$, 
$v^{\vee }$, $v'$ et $v'^{\vee }$.

\null
Preuve: \rm
(cf. [C] proposition 5.2.4)

\null
Le degr\'e formel de $\pi $ ne change pas si on remplace $\pi $ par un \'el\'ement de son orbite inertielle $\o $. On peut donc poser $d(\o ):=d(\pi )$.

\null
\null
\null
{\bf 2.} Fixons un couple parabolique standard $(P,M)$ de $G$ 
et une repr\'esentation cuspidale unitaire irr\'eductible $(\sigma ,E)$ de $M$. Notons $\o $ l'orbite inertielle 
de $\sigma $. Soit $\xi :\o \rightarrow i_{\overline {P}\cap K}^KE\otimes i_{P\cap K}^KE^{\vee }$ une application polynomiale \`a image dans un espace de dimension finie. Posons 
$\zeta (\sigma ')=(J_{\overline {P}\vert P}(\sigma ')^{-1}\otimes 1)\xi (\sigma ')$ 
pour tout $\sigma '\in \o$ et
$$f_{\zeta }(g)=\int _{\hbox {\srm Re}(\chi )=\mu >>_P0} 
E_{P,\sigma \otimes \chi }^G (\zeta (\sigma \otimes \chi ))(g^{-1}) d(Im(\chi )).$$

\null
{\bf 2.1} {\bf Proposition:} \it
La fonction $f_{\zeta }$ appartient \`a $C_c^{\infty }(G)$.

\null
Preuve: \rm
D'apr\`es 1.2.2 et 1.3, la fonction $f_{\zeta }$ est bien d\'efinie. Il est clair qu'elle est lisse. Il reste
donc \`a montrer que son support est compact.

Par la  d\'ecomposition de Cartan, on a
$$G=KM_0^+K \ \ \hbox{avec}\ \ M_0^+=\{m\in M_0\vert \ \forall \alpha \in \Delta \ <\alpha ,H_0(m)>\ \geq 0\}.$$
Par ailleurs, $KmK=Km'K$ si et seulement si $H_0(m)=H_0(m')$. Posons $A_0^+=A_0\cap M_0^+$. Il existe un sous-ensemble compact $C$ de $A_0^+$ tel que $M_0^+=CA_0^+$. Par un argument de $C$-lisset\'e,  il suffit de montrer que, pour tout $v\in i_{\overline {P}\cap K}^KE$, tout $v^{\vee }\in i_{P\cap K}^KE^{\vee }$  et toute fonction polynomiale $p$ sur $X (M)$, la fonction  sur $A_0^+$ d\'efinie par
$$a\mapsto \int _{{\srm Re}(\chi )=\mu >>_P0}\ p(\chi ) <i_P^G(\sigma \otimes \chi )(a)J_{\overline {P}\vert P}(\sigma \otimes \chi )^{-1}v,v^{\vee }> d(Im(\chi ))\eqno {\hbox {(*)}}$$
est \`a support compact.

Supposons d'abord $G$ semi-simple. Alors $\Delta $ est une base de $a_0^*$. Pour $\Theta \subseteq \Delta $ et $t', t\geq 0$, posons $A_0^+(\Theta ,t,t')=\{a\in A_0^+\vert <\alpha ,H_0(a)>\ \leq t\ \forall \alpha \in \Theta \  \hbox {\rm et}\  <\alpha ,H_0(a)>\ >t'\  \forall \alpha \in\Delta -\Theta \}$, et $A_0^+(\Theta ,t):=A_0^+(\Theta ,t,t)$.

On va montrer l'existence d'une fonction $(\Theta ,t)\mapsto f(\Theta ,t)$, $\Theta \subset \Delta $, $t\geq 0$, telle que la fonction (*) soit nulle en tout $a\in A_0^+(\Theta ,t,f(\Theta ,t))$. Ceci implique la proposition dans le
cas semi-simple:  Comme $A_0^+=\bigcup _{\Theta \subseteq \Delta }A_0^+(\Theta ,t)$ pour tout $t\geq 0$, il suffit d'en d\'eduire que, pour tout $\Theta \subseteq \Delta $ et tout $t\geq 0$, la restriction de (*) \`a $A_0^+(\Theta ,t)$ est \`a support compact. 

Effectuons une r\'ecurrence d\'ecroissante sur $\Theta \subseteq \Delta $. Comme $\Delta $ est une base de $a_0^*$, $A_0^+(\Delta ,t)$ est compact pour tout $t\geq 0$. Soit $\Theta \subset \Delta $. Si (*) est non nulle en $a\in A_0^+(\Theta ,t)$, alors $\Theta '=\{\alpha \in \Delta \vert <\alpha ,H_0(a) >\leq f(\Theta ,t)\}$ contient proprement $\Theta $. Par suite, $a \in \bigcup _{\Theta '\supset \Theta } A_0^+(\Theta ', f(\Theta ,t))$. Or, par hypoth\`ese de r\'ecurrence, la restriction de l'application (*) \`a cette r\'eunion est \`a support compact. 

Fixons $\Theta \subset \Delta $, $t\geq 0$ et montrons l'existence de $f(\Theta ,t)$. Posons $(P',M')=(P_{\Theta },M_{\Theta })$. Comme dans [W], on d\'eduit de la formule de Casselman (1.3.1), qu'il existe $t_0\geq t$, tel que, pour tout $a\in A_0^+(\Theta ,t,t_0)$, on ait
$$\eqalign {
&<i_P^G(\sigma \otimes \chi )(a)J_{\overline {P}\vert P}(\sigma \otimes \chi )^{-1}v, v^{\vee }>\cr
=&\ \gamma (G\vert M')^{-1}\delta _{P'}^{1/2}(a)\sum _{w\in W(M,M')} 
c_{P'\vert P}(\sigma \otimes \chi ,w)(J_{\overline {P}\vert P}(\sigma \otimes \chi )^{-1}v\otimes v^{\vee })(a).\cr }$$

En particulier, le coefficient matriciel est nul, si $W(M,M')=\emptyset$. 

L'\'etude de (*) se ram\`ene donc \`a celle de
$$a\mapsto \int _{\hbox {\srm Re}(\chi )=\mu >>_P0} p(\chi )\ c_{P'\vert P}(\sigma \otimes \chi ,w) (J_{\overline {P}\vert P}(\sigma \otimes \chi )^{-1}v\otimes v^{\vee })(a) d(Im(\chi ))\eqno {\hbox {(**)}}$$
pour tout $w\in W(M,M')$. Fixons $w\in W(M,M')$. 

Il existe un compact $C_{\Theta }\subset A_0$, $C_{\Theta }^{-1}\subset A_0^+$, tel que tout \'el\'ement de $A_0^+$ puisse s'\'ecrire sous la forme $a=a_{\Theta }a'c_a$ avec $a_{\Theta }\in A_{M'}$, $a'\in A_0\cap M'^{\hbox {\srm der}}$ et $c_a\in C_{\Theta }$. On a $<\alpha ,H_0(a')>\leq 0$ pour tout 
$\alpha \in \Delta -\Theta $, et l'ensemble des $H_0(a')$ avec $a\in A_0^+(\Theta ,t)$ est fini. En particulier, $a\in A_0^+(\Theta ,t,t')$ implique $a_{\Theta }\in A_{\Theta }\cap A_0^+(\Theta ,t,t')$.

Pour $a\in A_0^+$, notons $r_{a}$ la fonction rationnelle $\chi \mapsto p(\chi )c_{P'\vert P}(\sigma \otimes \chi ,w)(J_{\overline {P}\vert P}(\sigma \otimes \chi )^{-1}v\otimes v^{\vee })(a'c_a)$ d\'efinie sur $X (M)$. Par un argument de lisset\'e, on d\'eduit de la finitude de l'ensemble des $H_0(a')$, $a\in A_0^+(\Theta ,t)$, la finitude de l'ensemble des fonctions $r_{a}$, $a\in A_0^+(\Theta ,t)$.

Par ce qui pr\'ec\`ede et apr\`es avoir effectu\'e le changement de base $\chi \mapsto w^{-1}\chi $, l'\'etude de (**) se ram\`ene \`a celle de
$$a\mapsto \int _{\hbox {\srm Re}(\chi )=w\mu }\chi (a_{\Theta }) r_a(w^{-1}\chi ) d(Im(\chi )),\eqno {\hbox {(***)}}$$
o\`u $\mu $ a \'et\'e choisi suffisamment positif dans la chambre de Weyl de $P$ dans $a_M^*$.

On d\'eduit de 1.3 et de 1.3.2 que les p\^oles des fonctions rationnelles $\chi \mapsto r_a(w^{-1}\chi )$, $a\in A_0^+(\Theta ,t)$, (qui sont en nombre fini) sont de la forme $\chi _{\lambda }$ avec $\lambda $ sur un nombre fini de hyperplans de $a_M^*$ de la forme $<\alpha ^{\vee }, \lambda >=c$ avec $\alpha \in \Sigma (wPw^{-1})\cap \Sigma (wP_w'w^{-1})$.

A l'aide de la d\'ecomposition $a_{wMw^{-1}}^*=a_{M'}^*\oplus a_{wMw^{-1}}^{M' *}$, on d\'efinit $\chi _{\mu '}\in X (wMw^{-1})$ pour $\mu '\in a_{M'}^*$. Supposons $\mu '$ dans la chambre de Weyl positive de $P'$ dans $a_{M'}^*$. Soit $\alpha \in \Sigma (wPw^{-1})\cap \Sigma (wP_w'w^{-1})$. Alors, ou $\alpha \in \Sigma (wPw^{-1}\cap M')$ ou $\alpha _{\vert A_{M'}}\in \Sigma (P')$. Dans le premier cas, $<\mu ',\alpha ^{\vee }>=0$, alors que $<\mu ',\alpha ^{\vee }>\ >0$ dans le deuxi\`eme. L'int\'egrale dans (***) ne change donc pas de valeur pour $a\in A_0^+(\Theta ,t)$, si on remplace $w\mu $ par $w\mu +\mu '$ avec $\mu '$ dans la chambre de Weyl positive de $P'$ dans $a_{M'}^*$. 

L'ensemble des fonctions rationnelles $\chi \mapsto r_{a}(w^{-1}\chi )$,
$a\in A_0^+(\Theta ,t)$, \'etant fini, on peut choisir $t_w\geq t_0$ tel que, pour tout $a\in A_0^+(\Theta ,t,t_w)$,  $\chi \mapsto \chi (a_{\Theta })$ soit \`a d\'ecroissance rapide par rapport \`a $\chi \mapsto r_{a}(w^{-1}\chi )$,
lorsque Re$(\chi _{\vert A_{M'}})$ devient tr\`es positif dans la 
chambre de Weyl de $P'$ dans $a_{M'}^*$. 

Comme, pour $a\in A_0^+(\Theta ,t,t_w)$, l'int\'egrale dans (***) reste invariante si on remplace $w\mu $ par $w\mu +\mu '$, $\mu '>_{P'}0$, et que la fonction \`a l'int\'erieur de l'int\'egrale converge vers $0$ lorsque $\mu '$ devient tr\`es positif dans la chambre de Weyl de $P'$, (***) est nulle en tout $a\in A_0^+(\Theta ,t,t_w)$.

On pourra alors prendre pour $f(\Theta ,t)$ le plus grand des $t_w$.

\null
Consid\'erons maintenant le cas d'un groupe r\'eductif qui n'est pas semi-simple. On a $A_0=A_G(A_0\cap G^{\hbox {\srm der}})C'$ avec $C'$ compact. Les morphismes de restriction donnent lieu \`a une suite exacte $0\rightarrow X (G)\times X\rightarrow X (M) \rightarrow X (M\cap G^{\hbox
{\srm der}})\rightarrow 0$, o\`u $X$ d\'esigne un sous-ensemble fini de $X _{\hbox {\srm im}} (M)$ form\'e de caract\`eres de restriction triviale \`a $A_G (M\cap G^{\hbox{\srm der}})$. On identifie $X _{\hbox {\srm im}} (G)$ \`a un sous-groupe de $X _{\hbox {\srm im}} (M)$ au moyen de  ce morphisme. On a donc un isomorphisme $X _{\hbox {\srm im}} (M)/(X _{\hbox {\srm im}} (G)\times X)\rightarrow X _{\hbox {\srm im}} (M\cap G^{\hbox {\srm der}})$. 

Choisissons $\mu $ suffisamment positif dans la chambre de Weyl de $P$ dans $a_M^*$ et tel que $<\mu ,H>=0$ pour $H\in a_G$. Par ce qui pr\'ec\`ede et un argument de lisset\'e, l'\'etude de (*) se ram\`ene \`a celle de  
$$\eqalign {
\int _{\chi _{\mu } X_{\hbox {\ssrm im}} (M\cap G^{\hbox {\ssrm der}})}&\int _{X_{\hbox {\ssrm im}} (G)\times X}\ p(\chi '\chi _G)\chi _G(a_G)\cr 
&<i_P^G(\sigma \otimes \chi '\chi _G)(a')J_{\overline {P}\vert P}(\sigma
\otimes\chi '\chi _G)^{-1} v,v^{\vee }> d(Im(\chi _G)) d(Im(\chi '))\cr }
\eqno (\#)$$
pour $(a_G,a')$ dans $A_G\times (A_0^+\cap G^{\hbox{\srm der}})$. Il reste donc \`a montrer l'existence de compacts $C_G$ et $C_0$ de $A_G$ et $A_0^+\cap G^{\hbox{\srm der}}$ respectivement, tels que $(\#)$ soit nul si $(a_G,a')\not\in C_G\times C_0$. 

Remarquons d'abord que, si $a'\in A_0^+\cap G^{\hbox{\srm der}}$,
$$\eqalign {
&<i_P^G(\sigma \otimes \chi '\chi _G)(a')J_{\overline {P}\vert P}
(\sigma \otimes \chi '\chi _G)^{-1}v,v^{\vee }> \cr
=&<i_{P\cap G^{\hbox {\ssrm der}}}^{G^{\hbox {\ssrm der}}}
(\sigma \otimes \chi ')(a')
J_{\overline {P}\cap G^{\hbox {\ssrm der}}\vert P\cap G^{\hbox {\ssrm der}}}(\sigma \otimes \chi ')^{-1}v_{\vert G^{\hbox {\ssrm der}}}
,v^{\vee }_{\vert G^{\hbox {\ssrm der}}}>, \cr}$$
et que l'on peut remplacer ci-dessus $\int _{X_{\hbox {\ssrm im}} (G)\times X}$ par $\vert X \vert \ \int _{X_{\hbox {\ssrm im}} (G)}$. Comme $\chi '\mapsto \int _{X_{\hbox {\ssrm im}} (G)}\ p(\chi '$ $ \chi _G)$ $\chi _G(a_G)\ d(Im(\chi _G))$ est une application polynomiale sur $X (M\cap G^{\hbox {\srm der}})$ et que $M\cap G^{\hbox {\srm der}}$ est un sous-groupe de L\'evi semi-standard du groupe semi-simple $G^{\hbox {\srm der}}$, l'existence de $C_0$ r\'esulte du cas semi-simple consid\'er\'e pr\'ec\'edemment. L'int\'egrale sur $X _{\hbox {\srm im}} (G)$ portant sur une fonction polynomiale en $\chi _G$, l'existence de $C_G$ est imm\'ediate. \hfill {\fin 2}
  
\null
\null
{\bf 2.2} Lorsque $f$ est un \'el\'ement de $C_c^{\infty }(G)$ et que $P'=M'U'$ est un
sous-groupe parabolique de $G$, posons
$$f_{P'}(m')=\delta _{P'}(m')^{1/2}\int _{U'} f(m'u')du'\ \ \ \hbox{pour $m'\in M'$}.$$

\null
{\bf 2.2.1} {\bf Lemme:} \it 
On a $f_{P'}\in C_c^{\infty }(M')$. \rm

\null
{\bf 2.2.2} {\bf Lemme:} \it
Soit $(\pi ,V)$ une repr\'esentation lisse de $G$  et $v\in V$. Si $P'=M'U'$ est un sous-groupe parabolique et $H$ un sous-groupe ouvert compact de $G$ qui laisse $v$ invariant et qui admet une d\'ecomposition d'Iwahori $H=(U'\cap H)(M'\cap H)(\ol {U'}\cap H)$ par rapport au couple parabolique $(P',M')$, alors on a
$$\int _{U'\cap H}\pi (u'a)vdu'= {\mes (U'\cap H)\over \mes (H)}\int _H\pi (ha)v dh.$$
pour tout $a\in A_{M'}$ strictement positif pour $P'$. En particulier, l'\'el\'ement de $V$ \'egal \`a cette int\'egrale est invariant par $H$. 
\rm

\null
{\bf 2.2.3} {\bf Proposition:} \it
Soit $(P',M')$ un couple parabolique semi-standard, et supposons 

{\noindent $\dim M'\leq \dim M$.}
Alors on a

(i) $(f_{\zeta })_{P'}=0$ si $M'$ et $M$ ne sont pas conjugu\'es;

(ii) 
$$\eqalign {\gamma (G\vert M) &(f_{\zeta })_{P'}(m')=\sum _{w\in W(M,M)} \int _{\hbox {\srm Re}(w\chi )=\mu _w>>_{P'}0}\cr 
& E_{M,w(\sigma \otimes \chi )}^M
((\lambda (w)J_{P_w'\vert \overline {P}}(\sigma \otimes \chi )\otimes 
\lambda (w)J_{P_w'\vert P}(\sigma ^{\vee }\otimes \chi ^{-1}))\xi (\chi ))(m'^{-1})\ dIm(\chi )\cr
}$$ 
si $M'=M$.

\null
Preuve: \rm
Comme $\xi $ est une combinaison lin\'eaire d'applications de la forme $\chi \mapsto p(\chi )v\otimes v^{\vee }$, il suffit de montrer la proposition dans le cas o\`u $\xi $ est une telle fonction.

Soit $H$ un sous-groupe ouvert compact de $K$, admettant une d\'ecomposition d'Iwahori relative \`a $(P',M')$ et laissant $v$ et $v^{\vee }$ invariant. Soit $a\in A_{M'}$ strictement positif pour $P'$. Posons $U_0'=H\cap U'$. On a $U'=\bigcup _{l=0}^{\infty } a^{-l}U_0'a^l$. 
Par suite,
$$\eqalign {
&\delta _{P'}(m')^{-1/2}(f_{\zeta })_{P'}(m') \cr
=&\int _{U'}f_{\zeta }(m'u')du'\cr
=&\lim _{l\rightarrow \infty} \int _{a^{-l}U_0'a^l}\int _{\hbox {\srm Re}(\chi )=\mu
>>_P0}E_{P,\sigma \otimes \chi }^G(\zeta (\sigma \otimes \chi ))(u'^{-1}m'^{-1})dIm(\chi )\cr
=&\lim _{l\rightarrow \infty}\delta _{P'}(a^l)^{-1}\int _{U_0'}\int _{\hbox {\srm Re}(\chi
)=\mu >>_P0}p(\chi ) \cr
&<i_P^G(\sigma \otimes \chi )(a^l)J_{\overline {P}\vert P}(\sigma\otimes \chi )^{-1}
i_{\overline {P}}^G(\sigma \otimes \chi )(m'^{-1})v,i_P^G(\sigma^{\vee }\otimes\chi ^{-1})(u'a^l)v^{\vee }>dIm(\chi )du'.\cr  
}$$
En posant $v_{m'}=i_{\overline {P}}^G(\sigma \otimes \chi )(m'^{-1})v$ et 
$v_l^{\vee }=\int _{U_0'}i_P^G(\sigma ^{\vee }\otimes \chi ^{-1})(u'a^l)v^{\vee }du'$, ceci devient \`a l'aide du th\'eor\`eme de Fubini \'egal \`a
$$\lim _{l\rightarrow \infty }\delta _{P'}(a)^{-l}\int _{\hbox {\srm Re}(\chi)=\mu >>_P0}
p(\chi )<i_P^G(\sigma \otimes \chi )(a^l)J_{\overline {P}\vert P}(\sigma \otimes \chi
)^{-1}v_{m'},v_l^{\vee }>dIm(\chi ).$$

Il r\'esulte du lemme 2.2.2 que $v_l^{\vee }$ reste dans un espace de dimension finie pour $l>>0$. On peut
donc appliquer la formule de Casselman, et on trouve
$$\eqalign {
&<i_P^G(\sigma \otimes \chi )(a^l)J_{\overline {P}\vert P}(\sigma \otimes\chi
)^{-1}v_{m'},v_l^{\vee }>\cr
=&\gamma (G\vert M')^{-1}\delta _{P'}(a^l)^{1/2}\sum _{w\in W(M,M')} c_{P'\vert P}(\sigma \otimes
\chi ,w)(J_{\overline {P}\vert P}(\sigma \otimes\chi )^{-1}v_{m'}\otimes v_l^{\vee })(a^l)\cr
}$$ 
pour $l$ assez grand.

Le (i) de la proposition en r\'esulte aussit\^ot. Supposons dans la suite $M'=M$.

Observons que
$$\eqalign {
&(\lambda (w)J_{\ti {P}'_w\vert P}(\sigma ^{\vee }\otimes \chi ^{-1})v_l^{\vee })(1)\cr
=&\int _{U_0'}(J_{\ti {P}'_w\vert P}(\sigma ^{\vee }\otimes \chi ^{-1})i_P^G(\sigma ^{\vee }\otimes \chi ^{-1})(u'a^l)v^{\vee })(w^{-1})du'\cr
=&\int _{U_0'}(J_{\ti {P}'_w\vert P}(\sigma ^{\vee }\otimes \chi ^{-1})v^{\vee })(w^{-1}u'a^l)du'\cr
=&\delta _{P'}(a^l)^{1/2} w(\sigma ^{\vee }\otimes \chi ^{-1})(a^l)\int _{a^{-l}U_0'a^l}(J_{\ti
{P}'_w\vert P}(\sigma ^{\vee }\otimes \chi ^{-1})v^{\vee })(w^{-1}u')du'.\cr
}$$
Par suite,
$$\eqalign {
&\delta _{P'}(m')^{-1/2}(f_{\zeta })_{P'}(m')\cr
=&\gamma (G\vert M')^{-1}\lim _{l\rightarrow \infty }\sum _{w\in W(M,M)}\int _{\hbox {\srm Re}(\chi )=\mu >>_P0}p(\chi )\cr
&<(\lambda (w)J_{P_w'\vert P}(\sigma \otimes \chi)J_{\overline {P}\vert P}(\sigma \otimes
\chi)^{-1}v_{m'})(1),\cr
&\ \ \ \ \ \ \ \ \ \ \ \ \ \ \ \ \ \ \ \ \ \ \ \ \ \ \ \ \ \ \ \ \ \ 
\ \ \ \ \ \ \ \ \ \ \ 
\int _{a^{-l}U_0'a^l}(J_{\ti{P}'_w\vert P}
(\sigma ^{\vee }\otimes \chi ^{-1})v^{\vee })(w^{-1}u')du'>dIm(\chi ).\cr 
}$$

Fixons $w\in W(M,M)$ et calculons la limite correspondante. Il r\'esulte de 1.3 et de 1.3.2 qu'il existe un nombre fini de hyperplans de la forme $<\lambda ,\alpha ^{\vee }>=c$, $\alpha \in \Sigma (P)\cap \Sigma (P_w')$, tels que tout p\^ole de la fonction dans l'int\'egrale soit de la forme $\chi _{\lambda }$ avec $\lambda $ sur un de ces hyperplans. On peut donc remplacer ci-dessus $\mu >>_P0$ par $\mu _w>>_{P_w'}0$. Or, alors $\lim _{l\rightarrow \infty }\int_{a^{-l}U_0'a^l}(J_{\ti{P}'_w\vert P}(\sigma ^{\vee }\otimes \chi ^{-1})v^{\vee })(w^{-1}u')du'=(\lambda(w)J_{P_w'\vert \ti {P}_w'}(\sigma ^{\vee }\otimes \chi ^{-1})J_{\ti{P}'_w\vert P}(\sigma ^{\vee }\otimes \chi ^{-1})v^{\vee })(1)$. En appliquant 1.3.2 et la formule de produit pour les op\'erateurs d'entrelace-ment, on en d\'eduit que le terme correspondant \`a $w$ dans la somme ci-dessus est \'egal \`a
$$\int _{\hbox {\srm Re}(\chi )=\mu _w>>_{P_w'}0}p(\chi )<(\lambda (w)J_{P_w'\vert \overline{P}}(\sigma \otimes \chi )v_{m'})(1),(\lambda (w)J_{P_w'\vert P}(\sigma ^{\vee }\otimes \chi ^{-1})v^{\vee })(1)>dIm(\chi ),$$
d'o\`u la formule (ii). \hfill {\fin 2}

\null
\null
{\bf 2.3} {\bf Proposition:} \it
Soit $(\sigma ',E')$ une repr\'esentation irr\'eductible cuspidale d'un sous-groupe de L\'evi semi-standard
$M'$ de dimension inf\'erieure ou \'egale \`a celle de $M$. Soit $P'\in \hbox {\main P}(M')$.
Alors on a

(i) $\widehat {f_{\zeta }}(\sigma ',P')=0$ si $M'$ et $M$ ne sont pas conjugu\'es;

(ii) 
$$\eqalign {
&\gamma (G\vert M)\widehat {f_{\zeta }}(\sigma ',P')\cr
=&\ d(\sigma ')^{-1}\sum _{w\in W(M,M),\sigma '\in w\so } 
(J_{P'\vert \overline {wP}}(\sigma ')\lambda (w)\otimes J_{P'\vert wP}(\sigma '^{\vee })\lambda (w))\xi (w^{-1}\sigma ')\cr
}$$
si $M'=M$. \rm

\null
Pour la preuve de cette proposition, on utilisera les deux lemmes suivants, o\`u  $\lambda $ et $\rho $ d\'esignent respectivement l'action par translations \`a gauche et \`a droite de $G$ sur $C_c^{\infty }(G)$.

\null
{\bf 2.3.1} {\bf Lemme:} \it
Soit $(\sigma ',E')$ une repr\'esentation cuspidale de $M'$, $P'\in
\hbox {\main P}(M')$. Alors on a $\widehat {f}_{\zeta }^G(\sigma ',P')(g,h)=(\lambda (g)\rho (h)f_{\zeta })\widehat {_{P'}^{\ }}^{M'}(\sigma ',M')$.

\null
Preuve: \rm
Il suffit de combiner les lemmes VII.1.2 et VII.1.1 de [W].\hfill {\fin 2}

\null
\null
{\bf 2.3.2} {\bf Lemme:} \it 
On a $(\lambda (g)\rho (h)f_{\zeta })=f_{\zeta '}$, o\`u $\zeta
'(\chi )=(i_P^G(\sigma \otimes \chi )\otimes i_P^G(\sigma ^{\vee }\otimes \chi ^{-1}))(g,h)$ $\zeta (\chi )$.
\rm

\null
\it Preuve: \rm (de la proposition) 
Comme $\xi $ est une combinaison lin\'eaire d'applications de la forme $\chi \mapsto p(\sigma \otimes \chi )\ v\otimes v^{\vee }$ avec $p$ fonction polynomiale sur $\o $, il suffit de montrer la proposition dans le cas o\`u $\xi $ est une telle fonction.

La partie (i) de la proposition est une cons\'equence imm\'ediate des lemmes 2.2.3 (i) et 2.3.1.

Supposons $M'=M$. A l'aide de 2.2.3 (ii), on trouve avec $e'\in E'$ et $e'^{\vee }\in E'^{\vee }$,

$$\eqalign {
&\gamma (G\vert M)<(f_{\zeta })\widehat {_{P'}}^{M}(\sigma ',M)e',e'^{\vee }>\cr
=&\sum _{w\in W(M,M)} \int _M \int _{\hbox {\srm Re}(w\chi )=\mu _w>>_{P'}0}p(\sigma \otimes \chi )\cr
&<(\lambda (w)J_{P_w'\vert \overline {P}}(\sigma \otimes\chi )v)(1),w(\sigma \otimes \chi
)^{\vee }(m)(\lambda (w)J_{P_w'\vert P}(\sigma ^{\vee }\otimes \chi ^{-1})v^{\vee })(1)>dIm(\chi )\cr
&<\sigma '(m)e',e'^{\vee }>dm\cr 
=&\sum _{w\in W(M,M)}\int _{A_M\backslash M}<\sigma '(m)e',e'^{\vee }>\int _{A_M\cap K\backslash A_M}\int _{\hbox{\srm Re}(w\chi )=-\mu _w<<_{P'}0} p(\sigma \otimes \chi ^{-1})\cr
&<(\lambda (w)J_{P_w'\vert \overline {P}}(\sigma \otimes\chi ^{-1})v)(1),w(\sigma \otimes \chi ^{-1}
)^{\vee }(m)(\lambda (w)J_{P_w'\vert P}(\sigma ^{\vee }\otimes \chi )v^{\vee })(1)>\cr
&\ \ \ \ \ \ \ \ \ \ \ \ \ \ \ \ \ \ \ \ \chi _{\sigma '}(a)w(\chi \chi _{\sigma }^{-1})(a)(\int _{A_M\cap K}(\chi
_{\sigma '}(w\chi _{\sigma })^{-1})(\epsilon _a)d\epsilon _a)\ dIm(\chi )\ da\ dm.\cr 
}$$
Comme $A_M\cap K$ est compact de mesure 1, l'int\'egrale sur $A_M\cap K$ n'est non nulle que si ${\chi _{\sigma '}}_{\vert A_M\cap K}=w{\chi _{\sigma }}_{\vert A_M\cap K}$, et sa valeur est alors $1$. Dans ce cas, $\chi _{\sigma '}(w\chi _{\sigma })^{-1}$ est la restriction \`a $A_M$ d'un certain \'el\'ement $w\chi _w$ de $X (M)$.

A l'aide de la th\'eorie de Fourier sur un tore (cf. 1.2.3), on trouve alors 

$$\eqalign {
=&\sum _{w\in W(M,M)}\sum _{\chi \in \ker (X_{\hbox {\ssrm im}} (M)\rightarrow X_{\hbox {\ssrm im}} (A_M))} p(\sigma \otimes \chi _w\chi )\int _{A_M\backslash
M}<\sigma '(m')e',e'^{\vee }>\cr &<(\lambda (w)J_{P_w'\vert \overline {P}}(\sigma \otimes\chi
_w\chi )v)(1), \cr
&\ \ \ \ \ \ \ \ \ \ w(\sigma \otimes \chi _w \chi )^{\vee }(m)(\lambda (w)J_{P_w'\vert P}(\sigma ^{\vee }\otimes \chi _w^{-1}\chi ^{-1})v^{\vee })(1)>dm\cr 
=\ &d(\sigma
')^{-1}\sum _{w\in W(M,M), w^{-1}\sigma '\in \so}p(w^{-1}\sigma ')<(\lambda (w)J_{P_w'\vert 
\overline {P}}(w^{-1}\sigma ')v)(1),e'^{\vee }>\cr
&\ \ \ \ \ \ \ \ \ \ \ \ \ \ \ \ \ \ \ \ \ \ \ \ \ \ \ \ \ \ \ \ \ \ \ \ \ \ \ \ \ \ \ \ \ \ \ \ \ \ \ \ \ \ \ \ \ \ \ \ \ \ \ \  
<e',(\lambda(w)J_{P_w'\vert P}(w^{-1}\sigma '^{\vee })v^{\vee })(1)>\cr
}$$
par 1.4, si Re$(w\chi _w)<_{P'}-\mu _w$ pour tout $w$. Comme les deux applications sont rationnelles sur l'orbite inertielle de $\sigma '$, on a l'\'egalit\'e partout. 

A l'aide des lemmes 2.2.1, 2.2.3, 2.3.1, et 2.3.2, on en d\'eduit le r\'esultat \'enonc\'e.
\hfill {\fin 2}

\null
\null
{\bf 2.4} {\bf Proposition:} \it
Soit $\o '$ une orbite inertielle d\'efinie relative \`a un sous-groupe de L\'evi semi-standard 
de $G$. Supposons que $\o '$ et $\o $ ne soient pas
conjugu\'ees. Soit $\zeta '$ une application rationnelle sur $\o '$ v\'erifiant les propri\'et\'es
analogues \`a celles de $\zeta $ relatives \`a $\o '$. 

Alors on a
$$\int _G f_{\zeta }(g)\ \overline {f_{\zeta '}(g)}\ dg=0.$$

\null
Preuve: \rm
Les arguments dans la d\'emonstration de la proposition VII.2.2 dans [W] se g\'en\'eralisent, apr\`es
avoir remarqu\'e que pour tout $v\in i_{P\cap K}^KE$ et tout $v^{\vee }\in i_{P\cap K}^KE^{\vee }$ il
existe $v_1\in i_{P\cap K}^KE$ et $v_1^{\vee }\in i_{P\cap K}^KE^{\vee }$, tels que 
$\overline {<i_P^G(\sigma \otimes \chi )(g)v,v^{\vee }>}=
<i_P^G(\sigma \otimes \overline {\chi }^{-1})(g^{-1})v_1, v_1^{\vee }>$
pour tout $\chi \in X (M)$ et tout $g\in G$. 

\hfill {\fin 2}

\null
\null
{\bf 2.4.1} {\bf Corollaire:} \it
Soit $(\sigma ',E')$ une repr\'esentation irr\'eductible cuspidale d'un sous-groupe de L\'evi
$M'$ de $G$. Supposons que $\sigma '$ ne soit conjugu\'ee \`a aucun \'el\'ement de $\o $. 

Alors on a $\widehat {f_{\zeta }}(\sigma ',P')=0$.

\null
Preuve: \rm (cf. Corollaire VII.2.3 dans [W].)\hfill {\fin 2}

\null
\null
{\bf 3.} Notons $\Theta $ l'ensemble des couples $(P, \o)$ form\'es d'un sous-groupe parabolique semi-standard $P=MU$ et de l'orbite inertielle d'une repr\'esentation irr\'eductible cuspidale de $M$. Notons $\H (\Theta )$ l'ensemble des familles $\varphi =\{\varphi _{P, \so}\}_{(P,\so )\in \Theta }$ dont les composantes sont des applications polynomiales qui v\'erifient les conditions 1) - 4) du th\'eor\`eme 0.1.

\null
{\bf 3.1} Soit $\varphi \in \H(\Theta )$. Choisissons pour tout couple $(P,\o )$ une application polynomiale $\xi _{P,\so }$ v\'erifiant les conclusions de la proposition 0.2 relatives \`a $\varphi _{P,\so }$. Notons $\zeta _{P,\so }$ l'application rationnelle donn\'ee par $\zeta _{P,\so }(\sigma )=(J_{\overline {P}\vert P}(\sigma )^{-1}\otimes 1)\xi _{P,\so }(\sigma )$ pour $\sigma \in \o $.

\null
{\bf 3.1.1} {\bf Proposition:} \it
Soient $(P,\o )$, $(P',\o ')$ dans $\Theta $ et $\sigma '\in \o '$.

Alors 
$$\widehat {f_{\zeta _{P,\sso }}}(P',\sigma ')=\gamma (G\vert M')^{-1}d(\sigma ')^{-1}
\cases {\varphi _{P',\so '}(\sigma '), &si $\o '$ et $\o $ sont conjugu\'es;\cr
        0 ,&sinon. \cr}$$
       
\null
Preuve: \rm
Le corollaire 2.4.1  prouve que $\widehat {f_{\zeta _{P,\sso }}}(P',\sigma ')=0$ si $\o '$ et $\o $ ne sont pas conjugu\'es.

Si $\sigma '\in \o $ et $P=P'$, l'\'egalit\'e $\widehat {f_{\zeta _{P,\sso }}}(P',\sigma ')=\gamma (G\vert M')^{-1}d(\sigma ')^{-1}\varphi _{P',\so '}(\sigma ')$ r\'esulte imm\'ediatement des propositions 2.3 (ii) et 0.2.

Si $P\ne P'$, celle-ci se d\'eduit de l'identit\'e $(J_{P'\vert P}(\sigma ')\otimes J_{P\vert P'}(\sigma '^{\vee })^{-1})\varphi _{P,\sigma '}=\varphi
_{P',\sigma '}.$

Si finalement $\sigma '\in w\o $, $w\in W^G$, l'identit\'e $(\lambda (w)\otimes \lambda (w))\widehat {f_{\zeta _{P,\sso }}}(w^{-1}P',w^{-1}\sigma ')=\widehat {f_{\zeta _{P,\sso }}}(P',\sigma ')$ et ce que l'on vient de prouver dans le cas $\sigma '\in \o$
impliquent que $\widehat {f_{\zeta _{P,\sso }}}(P',\sigma ')$ $=\varphi _{P',\so '}(\sigma ')$. \hfill {\fin 2}

\null
\null
{\bf 3.1.2} {\bf Corollaire:} \it
La fonction $f_{\zeta _{P,\sso }}$ ne d\'epend pas du choix de $\xi _{P,\so }$.

\null
Preuve: \rm
Ceci r\'esulte de la proposition 3.1.1, puisque l'ensemble des transform\'ees de Fourier $\{f_{\zeta _{P,\sso }}(P',\sigma ')\}_{(P',\sigma ')}$ d\'etermine la fonction $f_{\zeta _{P,\sso }}$ (cf. [BZ] proposition 2.12).  \hfill {\fin 2}

\null
On pourra donc \'ecrire $f_{\varphi _{P,\sso }}$ \`a la place de $f_{\zeta _{P,\sso }}$.

\null
\null
{\bf 3.1.3} {\bf Corollaire:} \it
L'\'egalit\'e $f_{\varphi _{P,\sso }}=f_{\varphi _{P',\sso '}}$ vaut si $\o $ et $\o '$ sont conjugu\'es. 

\null
Preuve: \rm
La preuve est analogue \`a celle du corollaire 3.1.2. \hfill {\fin 2}

\null
\null
{\bf 3.2} Lorsque $(P=MU, \o )\in \Theta $, posons $[\o ]=\{w\sigma \vert w\in W^G, \sigma \in \o \}$, $W(M, \o )=\{ w\in W(M,M)\vert w\o =\o \}$, $c([\o ])= \vert \hbox {\main P} (M)\vert ^{-1} \vert W^M \vert \vert W^G \vert ^{-1}\vert W(M, \o )\vert\ \gamma (G\vert M)d(\o )$ et $f_{\varphi _{[\sso ]}}=c([\o ])\sum _{(P',\so ')\in \Theta ,\ \so '\subseteq [\so ]} f_{\varphi _{P',\sso '}}$.

Il r\'esulte de 3.1.1 que la fonction $f_{\varphi _{[\sso ]}}$ v\'erifie les propri\'et\'es 3) et 4) \'enonc\'ees dans l'introduction. La propri\'et\'e 5) est une cons\'equence directe de 2.4.

Apr\`es avoir rappel\'e qu'il est montr\'e dans [BZ] qu'un \'el\'ement de $C_c^{\infty }(G)$ est d\'etermin\'e par ses transform\'ees de Fourier, on s'aper\c coit que l'on a montr\'e le r\'esultat suivant:

\null
{\bf Th\'eor\`eme:} \it
Soit $\varphi $ dans $\H (\Theta )$. La fonction
$$f_{\varphi }=\sum _{(P,\so )\in \Theta } c([\o ]) f_{\varphi _{P, \sso }},$$
est l'unique \'el\'ement de $C_c^{\infty }(G)$ qui v\'erifie $\widehat {f}_{\varphi }(P, \sigma )=\varphi _{P, \so }(\sigma )$ pour tout $(P, \o )\in \Theta $, $\sigma \in \o $. \rm

\null
\null
{\bf 3.2.2} {\bf Corollaire:} \it Pour que $\varphi $ soit un \'el\'ement de $\H (\Theta )$, il faut et il suffit qu'il existe $f$ dans $C_c^{\infty }(G)$, telle que $\varphi _{P, \so }(\sigma )=\widehat {f}(P, \sigma )$ pour tout $(P, \o )\in \Theta $, $\sigma \in \o $. \rm

\null
\null
\null
{\bf B. Une relation polynomiale}

\vskip 1cm
{\bf 1.1} Si $(\pi ,V)$ est une repr\'esentation lisse de $G$ et $P=MU$ un sous-groupe parabolique de $G$, notons $\pi _P$ la repr\'esentation lisse de $M$ dans le module de Jacquet $V_P$ de $V$. Lorsque $U_1$ est un sous-groupe ouvert compact de $U$, \'ecrivons $V(U_1)$ pour l'ensemble des \'el\'ements $v\in V$ tels que $\int _{U_1}\pi (u)v du=0$. Le noyau $V(U)$ de la projection canonique $j_P:V\rightarrow V_P$ est la r\'eunion des $V(U_1)$, $U_1$ parcourant les sous-groupes ouverts compacts de $U$.

L'ensemble des \'el\'ements de $V$ invariants pour l'action par un sous-groupe ouvert $H$ de $G$ sera not\'e $V^H$.

\null
\null
{\bf 1.2} Soit $\o $ l'orbite inertielle d'une repr\'esentation irr\'eductible cuspidale de $M$. Rappelons que les fonctions rationnelles $j_{{\overline P}\vert P}$, $P\in \P(M)$, sont toutes \'egales \`a une m\^eme fonction, not\'ee $j$, et que tout point $W(M,M)$-r\'egulier de $\o $ est r\'egulier pour $j$.

Un \'el\'ement $\sigma $ de $\o $ sera dit \it en position g\'en\'erale \rm s'il v\'erifie les deux propri\'et\'es suivantes:

(i) le caract\`ere central de $\sigma $ est $W(M,M)$-r\'egulier;

(ii) $j(\sigma )\ne 0$.

\null
Fixons $\sigma \in \o$. L'ensemble des $\chi \in X (M)$ avec $\sigma \otimes \chi $ en position g\'en\'erale est Zariski dense dans $X (M)$. Deux applications rationnelles sur $\o $ sont donc \'egales d\`es qu'elles co\"\i ncident sur l'ensemble des points en position g\'en\'erale.

\null
\null
{\bf 1.3 Proposition:} \it
Soient $(\sigma ,E)$ une repr\'esentation irr\'eductible cuspidale de $M$ et $P'\in \P(M)$. Supposons $\sigma $ en position g\'en\'erale. Alors l'application
$$i_{P'}^GE\longrightarrow \bigoplus _{w\in W(M,M)}wE,\ \ \ \ \ 
v\mapsto\bigoplus _{w\in W(M,M)}\ ((J_{P\vert wP'}(w\sigma )\lambda (w))v)(1)$$
se factorise par $(i_{P'}^GE)_P$. Elle induit un isomorphisme
$$(i_{P'}^G\sigma )_P\longrightarrow \bigoplus _{w\in W(M,M)}\ w\sigma .$$

\null
Preuve: \rm Le r\'esultat \'enonc\'e relatif au module de Jacquet faible dans [W] au cours de la preuve de V.1.1 se g\'en\'eralise sans probl\`eme.
\hfill {\fin 2}

\null
\null
{\bf 2.} Soit $B$ l'anneau des fonctions r\'eguli\`eres d'une vari\'et\'e alg\'ebrique affine complexe. La notion d'une $B$-famille de repr\'esentations admissibles a \'et\'e d\'efinie dans [BD]. Le r\'esultat suivant a \'et\'e montr\'e par Casselman [C] dans le cas d'une repr\'esentation admissible. Sa preuve se g\'en\'eralise au cas d'une $B$-famille de repr\'esentations admissibles.

\null
{\bf Proposition:} \it
Soient $(\pi _B, V_B)$ une $B$-famille de repr\'esentations admissibles de $G$, $P=MU$ un sous-groupe parabolique semi-standard, et $H$ un sous-groupe ouvert compact de $G$ admettant une d\'ecomposition d'Iwahori par rapport \`a $(P,M)$.

Alors il existe un sous-groupe ouvert compact $U_1$ de $U$ tel que $V_B^H\cap V_B(U)\subseteq V_B(U_1)$. Les espaces $(V_B^H)_a:=\pi _B(1_{HaH})V_B$ avec $a\in A_M$ positif pour $P$ et v\'erifiant $aU_1a^{-1}\subseteq H\cap U$ sont tous \'egaux \`a un m\^eme espace, not\'e $S_P^H(V_B)$. Le foncteur de Jacquet induit un isomorphisme $S_P^H(V_B)\rightarrow (V_B)_P^{H\cap M}$ de $B$-modules.
\rm

\null
\null
{\bf 3.} Fixons un couple parabolique semi-standard $(P,M)$ et une repr\'esentation irr\'educti-ble cuspidale $(\sigma ,E)$ de $M$. Choisissons un sous-groupe ouvert compact $H$ de $G$ admettant une d\'ecomposition d'Iwahori par rapport \`a tout couple parabolique semi-standard. (On peut en trouver aussi petit que l'on veut.)

Notons $B=B_M$ l'anneau des fonctions r\'eguli\`eres sur la vari\'et\'e alg\'ebrique $X(M)$. Comme dans [W], on d\'eduit de $(\sigma ,E)$ les $B$-familles de repr\'esentations admissibles $(\sigma _B,$ $E_B)$ et $(\pi _B,$ $V_B)=(i_P^G\sigma _B, i_P^GE_B)$ de $M$ et $G$ respectivement. La classe d'isomorphie de $(\sigma _B,E_B)$ et $(\pi _B,V_B)$ ne change pas si on remplace $\sigma $ par un \'el\'ement de sa classe inertielle. On pourra donc \'ecrire $i_P^GE_{\so, B}$, si on ne veut distinguer aucun \'el\'ement de $\o $.

\null
\null
{\bf 3.1} Pour $\chi \in X(M)$, notons $E_{\chi }$ et $V_{\chi }$ respectivement les espaces des repr\'esentations $\sigma \otimes \chi $ et $i_P^G(\sigma \otimes \chi )$. On dispose de morphismes de sp\'ecialisation $\hbox{\rm sp}_{\chi }:E_B\rightarrow E_{\chi }$ et $\hbox{\rm sp}_{\chi }:V_B\rightarrow V_{\chi }$ qui commutent avec l'action du groupe. Ainsi, toute application polynomiale sur $X (M)$ \`a valeurs dans $E$ ou $i_{P\cap K}^KE$ correspond \`a un \'el\'ement de $E_B$ ou $V_B$, et vice versa. On \'ecrira \'egalement $E_{\sigma '}$, $V_{\sigma '}$ et $\sp _{\sigma '}$ si on ne veut distinguer aucun \'el\'ement de $\o $.

Le lemme suivant est une cons\'equence imm\'ediate des d\'efinitions:

\null
\null
{\bf Lemme:} \it 
Soit $P'\in \P(M)$. L'\'egalit\'e $\hbox {\rm sp}_{\chi }(S_{P'}^H(V_B))=S_{P'}^H(V_{\chi })$ vaut pour tout $\chi \in X(M)$. \rm

\null
\null
{\bf 3.2} Pour $b\in B$, notons $b^{\vee }$ l'\'el\'ement de $B$ qui v\'erifie $b^{\vee }(\chi )=b(\chi ^{-1})$. On d\'esignera par $E^{\vee }_{B^{\vee }}$ l'espace $E^{\vee }_B$ muni de la structure de $B$-module pour laquelle la multiplication scalaire $B\times E_B^{\vee }\rightarrow E_B^{\vee }$ est donn\'ee par $(b,e_B^{\vee })\mapsto b^{\vee }e_B^{\vee }$. Le produit de dualit\'e $<,>$ sur $E\times E^{\vee }$ induit par extension des scalaires une forme $B$-lin\'eaire $M$-\'equivariante $<,>_B$ sur $E_B\times E^{\vee }_{B^{\vee }}$. On en d\'eduit une forme bilin\'eaire $G$-\'equivariante $<,>_B$ sur $V_B\times V^{\vee }_{B^{\vee }}$. Pour $\chi \in X (M)$, $\sp _{\chi }(<,>_B)$ induit alors par passage au quotient le produit de dualit\'e entre $V_{\chi }$ et $V^{\vee }_{\chi ^{-1}}$.

\null
{\bf Proposition:} \it
Soit $P'\in \P(M)$. Les $B$-modules $S_{P'}^H(V_B)$ et $S_{\overline {P'}}^H(V^{\vee }_{B^{\vee }})$ sont libres de type fini et en dualit\'e par $<,>_B$.

\null
Preuve: \rm
En compl\'etant une base de $E^{H\cap M}$ en une base de $E$, on voit que $E_B^{H\cap M}$ est un $B$-module libre de rang \'egal \`a la dimension de $E^{H\cap M}$. On sait que $(V_B)_{P'}^{H\cap M}$ poss\`ede une filtration finie dont les sous-quotients sont des $B$-modules libres isomorphes \`a $E_B^{H\cap M}$. On en d\'eduit que $(V_B)_{P'}^{H\cap M}$ est libre de m\^eme rang que les espaces $(V_{\chi })_{P'}^{H\cap M}$ ou $(V_{\chi ^{-1}}^{\vee })_{P'}^{H\cap M}$, $\chi \in X (M)$. Il en est de m\^eme pour $(V_B)_{\overline {P'}}^{H\cap M}$.

Fixons des bases $\{ v_i\} _{i\in I}$ et $\{ v_i^{\vee }\}_{i\in I}$ de $S_{P'}^H(V_B)$ et $S_{\overline {P'}}^H(V_{B^{\vee }}^{\vee })$. Il suffit de montrer que la matrice $(<v_i, v_j^{\vee }>_B)_{i,j}$ est inversible, i.e. que son d\'eterminant $d$ appartient \`a $B^{\times }$. Or, dans le cas contraire $d$ serait contenu dans un id\'eal maximal $m_{\chi }$ de $B$ correspondant \`a un point $\chi $ de $X (M)$. En sp\'ecialisant en $\chi $, il en r\'esulterait que la forme bilin\'eaire $\sp _{\chi }(<,>_B)$ restreinte \`a $S_{P'}^H(V_{\chi })\times S_{\overline {P'}}^H(V_{\chi ^{-1}}^{\vee })$ serait d\'eg\'en\'er\'ee. Ceci est faux (cf. [C] th\'eor\`eme 4.2.4). \hfill {\fin 2}

\null
\null
{\bf 3.3} Notons $\o $ l'orbite inertielle de $\sigma $. Rappelons que toute repr\'esentation lisse $E'$ de $M$ admet une d\'ecomposition $E'=E'(\o )\oplus E'(\hors \o)$, telle que tout sous-quotient de $E'(\o )$ soit dans $\o $ et qu'aucun sous-quotient de $E'(\hors \o)$ ne le soit.

Pour $P'\in \P(M)$, notons $S_{P'}^H(V_B)(\o )$ le sous-$B$-module de $S_{P'}^H(V_B)$, form\'e des \'el\'e-ments \`a image dans $(V_B)_{P'}^{H\cap M}(\o )$. D\'efinissons de fa\c con analogue $S_{\overline {P'}}^H(V^{\vee }_{B^{\vee }})(\o ^{\vee })$,\ $S_{P'}^H(V_{\chi })$ $(\hors \o )$ etc.

\null
\null
{\bf 3.3.1 Lemme:} \it
Les espaces $S_{P'}^H(V_B)(\o )$ et $S_{\overline {P'}}^H(V^{\vee }_{B^{\vee }})(\hors \o^{\vee })$ sont orthogonaux.

\null
Preuve: \rm
Soient $v\in S_{P'}^H(V_B)(\o )$ et $v^{\vee }\in S_{\overline {P'}}^H(V^{\vee }_{B^{\vee }})(\hors \o^{\vee })$. Il suffit de montrer que $\sp _{\sigma ' }(<v,v^{\vee }>_B)=0$, si $\sigma '\in \o $ est en position g\'en\'erale. Fixons un tel $\sigma '$. On est donc ramen\'e \`a montrer que $<v,v^{\vee }>=0$ pour $v\in S_{P'}^H(V_{\sigma '})(\o )$ et $v^{\vee }\in S_{\overline {P'}}^H(V^{\vee }_{\sigma '^{\vee }})(\hors \o^{\vee })$. Rappelons (cf. [C] paragraphe 4) qu'il existe un produit bilin\'eaire $M$-\'equivariant $<,>_{P'}$ sur $j_{P'}(V_{\sigma '})\times j_{\ol{P'}}(V^{\vee }_{\sigma '^{\vee }})$ tel que 
$$<v,v^{\vee }>=<j_{P'}(v),j_{\ol {P'}}(v^{\vee })>_{P'}.$$
Comme $j_{P'}(V_{\sigma '})\simeq \bigoplus _{w\in W(M,M)} wE_{\sigma '}$ et que $j_{\ol {P'}}(V^{\vee }_{\sigma '^{\vee }})\simeq \bigoplus _{w\in W(M,M)} wE^{\vee }_{\sigma '^{\vee }}$ par choix de $\sigma '$, il suffit de consid\'erer le cas $j_{P'}(v)\in wE_{\sigma '}$ et $j_{\ol{P'}}(v^{\vee })\in w'E_{\sigma '^\vee }^{\vee }$. On a n\'ecessairement $w\ne w'$. Les espaces $wE_{\sigma }$ et $w'E_{\sigma '}$ n'ont donc pas d'entrelacement, d'o\`u $<j_{P'}(v),j_{\ol {P'}}(v^{\vee })$ $>_{P'}=0$. \hfill {\fin 2}

\null
\null
{\bf 3.3.2 Corollaire:} \it
Les $B$-modules libres $S_{P'}^H(V_B)(\o )$ et $S_{\overline {P'}}^H(V^{\vee }_{B^{\vee }})(\o ^{\vee })$ sont en dualit\'e par $<,>_B$.

\null
Preuve: \rm
Il suffit de rappeler que $S_{P'}^H(V_B)=S_{P'}^H(V_B)(\o )\oplus S_{P'}^H(V_B)(\hbox{\rm hors} \o)$. \hfill {\fin 2}

\null
\null
{\bf 3.4} On suppose dans cette section que $E^{H\cap M}\ne 0$.

\null
{\bf 3.4.1 Lemme:} \it
Soit $P'\in \P(M)$. Tout sous-quotient non nul $V'$ de $V_B$ en tant que $G$-module v\'erifie $V'_{P'}(\o )\ne 0$.

\null
Preuve: \rm
D'apr\`es [BD] proposition 2.8, on a une injection
$$V'\rightarrow \bigoplus _{P''\in \sP(M)} i_{P''}^G(V'_{P''}(\o )).$$
Par suite, $V'_{P''}(\o )\ne 0$ pour au moins un $P''\in \P (M)$. Ceci \'equivaut \`a dire qu'il existe $\sigma ''\in \o $ avec
$$0\ne \Hom _M(V'_{P''},\sigma '')=\Hom _G(V',i_{P''}^G\sigma '').$$
D\'eduisons-en $\sigma '\in \o$ avec $\Hom _M(V'_{P'},\sigma ')\ne 0$: Par it\'eration, on se ram\`ene \`a $P''$ et $P'$ adjacents. Il faut distinguer deux cas:

Si $wP''\ne P'$ ou $w\sigma ''\not\in \o$ pour tout $w\in W(M,M)$, l'op\'erateur d'entrelacement $J_{P'\vert P''}(\sigma '')$ est bien d\'efini et inversible, d'o\`u par composition un \'el\'ement non nul de $\Hom _G(V', i_{P'}^G\sigma '')=\Hom _M(V'_{P'},\sigma '')$.

Dans le cas contraire, choisissons $w\in W(M,M)$ et $\sigma '\in \o$ avec $i_{wP''}^G(w\sigma '')=i_{P'}^G\sigma '$. En composant avec $\lambda (w)$, on trouve donc un \'el\'ement non nul de $\Hom _G(V', i_{P'}^G\sigma ')$.

Ceci prouve bien que $V'_{P'}(\o )\ne 0$. \hfill {\fin 2}

\null
\null
{\bf 3.4.2 Proposition:} \it
Tout sous-$G$-module $V'$ de $V_B$ est engendr\'e par $V'\cap S_{P'}(V_B^H)(\o )$. 

\null
Preuve: \rm
Notons $V''$ le sous-$G$-module engendr\'e par cet ensemble. On a l'\'egalit\'e ${V''_{P'}}^{H\cap M}(\o )={V'_{P'}}^{H\cap M}(\o )$, d'o\`u $(V'/V'')_{P'}^{H\cap M}(\o )=0$ par exactitude du foncteur de Jacquet. Par choix de $H$, on d\'eduit de 3.4.1 que $V'=V''$. \hfill {\fin 2}

\null
\null
{\bf 4.} Rappelons que les ensembles $\Theta $ et $\H (\Theta )$ ont \'et\'e d\'efinis en A.3. On dira qu'un \'el\'ement $\varphi $ de $\H (\Theta )$ a la propri\'et\'e $(\P)$ en $(P, \o)$, si $\varphi _{P,\so }$ v\'erifie les conclusions de la proposition 0.2. On dira \'egalement que $\varphi _{P,\so }$ a la propri\'et\'e $(\P )$.

Fixons pour tout $(P,\o )\in \Theta $ un sous-groupe ouvert compact distingu\'e $H=H(\o )$ de $K$ admettant une d\'ecomposition d'Iwahori par rapport \`a tout couple parabolique semi-standard et tel que tout \'el\'ement de $\o $ admette des invariants par rapport \`a $H\cap M$. On peut par ailleurs supposer $H(\o ')=H$ si $\o $ et $\o '$ sont conjugu\'es, ce que l'on fera d\'esormais.

Pour $(P', \o ')\in \Theta $, $\o '$ conjugu\'e \`a $\o $, on peut donc parler, gr\^ace \`a 3.3.2, de la projection $\varphi _{P,\so }^{H,P',\so '}$ de $i_P^GE_{\so ,B}$ sur $S_{P'}^H(i_P^GE_{\so ,B})(\o ')$ de noyau \'egal \`a l'intersection des noyaux des \'el\'ements de $S_{\overline {P'}}^H(i_P^GE_{\so ^{\vee },B^{\vee }})(\o '^{\vee })$. Si $(P',\o')$ $=(P,\o)$, on \'ecrira plus simplement $\varphi _{P, \so }^H$.

\null
\null
{\bf 4.1} Soit $(P, \o)\in \Theta $, $H=H(\o )$. Rappelons que tout \'el\'ement $\varphi _B$ de $E_{\so ,B}\otimes _BE_{\so ^{\vee },B^{\vee }}$ ou de $i_P^GE_{\so ,B}\otimes _Bi_P^GE_{\so ^{\vee },B^{\vee }}$ correspond  par l'application de sp\'ecialisation $\sp _{\chi }$ \`a une application polynomiale sur $X (M)$ \`a valeurs respectivement dans $E_{\so }\otimes E_{\so ^{\vee }}$ ou dans $i_{P\cap K}^KE_{\so }\otimes i_{P\cap K}^KE_{\so ^{\vee }}$, et vice versa. On \'ecrira $\varphi _B(\sigma )=\sp _{\sigma }\varphi _B$. 

\null
{\bf 4.1.1 Lemme:} \it
Soient $P''\in \P(M)$, $\sigma \in \o $ et $w\in W^G$. Supposons $\sigma $ en position g\'en\'erale. Alors on trouve pour $(P',\o ')\in \Theta $, $\o '$ conjugu\'e \`a $\o $,

(i) $J_{P''\vert P}(\sigma )(S_{P'}^H(i_P^GE_{\sigma })(\o'))= S_{P'}^H(i_{P''}^GE_{\sigma })(\o ');$

(ii) $\lambda (w)(S_{P'}^H(i_P^GE_{\sigma })(\o '))=S_{P'}^H(i_{wP}^GwE_{\sigma })(\o ').$

\null
Preuve: \rm Il est clair par d\'efinition de $S_{P'}^H$ que l'image de $S_{P'}^H(i_P^GE_{\sigma })(\o')$ par les op\'erateurs $J_{P''\vert P}(\sigma )$ et $\lambda (w)$ est contenu dans $S_{P'}^H(i_{P''}^GE_{\sigma })$ et $S_{P'}^H(i_{wP}^GwE_{\sigma })$ respectivement. Il r\'esulte de la proposition 1.3 que tout sous-quotient de cette image est en fait un \'el\'ement de $\o '$. Ceci prouve l'inclusion. Les op\'erateurs consid\'er\'es \'etant bijectifs pour $\sigma $ en position g\'en\'erale, l'inclusion inverse s'en d\'eduit par un argument de sym\'etrie. \hfill {\fin 2}

\null
{\bf 4.1.2 Lemme:} \it Avec les notations du lemme {\rm 4.1.1}, on trouve

(i) $J_{P''\vert P}(\sigma )\varphi _{P,\so }^{H,P',\so '}(\sigma )=\varphi _{P'',\so }^{H,P',\so '}(\sigma )J_{P''\vert P}(\sigma )$;

(ii) $\lambda (w)\varphi _{P,\so }^{H,P',\so '}(\sigma )=\varphi _{wP,w\so }^{H,P',\so '}(w\sigma )\lambda (w).$

\null
Preuve: \rm Cela r\'esulte essentiellement du lemme pr\'ec\'edent et en appliquant la proposition 1.3. \hfill {\fin 2}

\null
\null
{\bf 4.2} {\bf Lemme:} \it
Pour tout $(P, \o)\in \Theta $, il existe un \'el\'ement $\varphi $ dans $\H (\Theta )$ avec $\varphi _{P,\so }=\varphi _{P,\so }^{H(\so )}$.

\null
Preuve: \rm
Si $\o '=w\o w^{-1}$ avec $w\in W$, posons $\varphi _{P',\so '}=\varphi _{P', \so '}^{H,P,\so }$. Sinon $\varphi _{P', \so '}=0$.

Les propri\'et\'es 1) et 2) du th\'eor\`eme 0.1 sont v\'erifi\'ees par d\'efinition de $\varphi _{P', \so '}^{H,P,\so }$.

Les propri\'et\'es 3) et 4) r\'esultent des \'egalit\'es (i) et (ii) du lemme 4.1.2. \hfill {\fin 2}

\null
\null
{\bf 4.3} Posons $W(M, \o )=\{ w\in W(M,M)\vert w\o =\o \}$.

\null
{\bf Lemme:} \it
Supposons: si on a un \'el\'ement $\varphi $ de $\H (\Theta )$ et $(P, \o)$ v\'erifiant $\varphi _{P, \so }=\varphi _{P, \so }^H$ avec $H=H(\o )$, alors cet \'el\'ement v\'erifie ($\P $) en $(P, \o )$.

Alors la propri\'et\'e $(\P )$ est v\'erifi\'ee pour tout $\varphi $ de $\H (\Theta )$ en tout $(P,\o )$.

\null
Preuve: \rm 
Fixons $(P, \o)$. Posons $V_B\times V^{\vee }_{B^{\vee }}=i_P^GE_{\so ,B}\times i_P^GE_{\so ^{\vee }, B^{\vee }}$. Montrons d'abord que tout \'el\'ement $\varphi $ de $\H (\Theta )$ qui v\'erifie $\varphi _{P, \so}\in S_P^H(V_B)(\o )\otimes _BS_{\overline {P}}^H(V^{\vee }_{B^{\vee }})(\o ^{\vee })$ a la propri\'et\'e $(\P )$ en $(P, \o )$:

En effet, de 4.2 et de nos hypoth\`eses, il r\'esulte l'existence d'une application polynomiale $\xi _{P, \so}$ sur $\o $ telle que, pour $\sigma \in \o $,
$$\varphi _{P, \so }^H(\sigma )=\sum _{w\in W(M,\so )} (J_{P\vert \overline {wP}}(\sigma )\lambda (w)\otimes J_{P\vert wP}(\sigma ^{\vee })\lambda (w)) \xi _{P, \so }(w^{-1}\sigma ).$$
On en d\'eduit que
$$\eqalign {
\varphi _{P,\so }(\sigma )=&\varphi _{P,\so }(\sigma )\varphi _{P, \so }^H(\sigma )\cr
=&\sum _{w\in W(M,\so )} (J_{P\vert \overline {wP}}(\sigma )\lambda (w)\otimes J_{P\vert wP}(\sigma ^{\vee })\lambda (w)) \varphi _{\overline {P}, \so}(w^{-1}\sigma )\xi _{P, \so }(w^{-1}\sigma ),\cr }$$
le produit dans la premi\`ere ligne d\'esignant la composition dans l'alg\`ebre $\End (i_P^GE_{\sigma })$. Comme $\sigma \mapsto \varphi _{\overline {P}, \so }(w^{-1}\sigma )$ est polynomiale, on a une relation du type voulu.

L'ensemble des \'el\'ements $\varphi _{P, \so }$ de $V_B\otimes V^{\vee }_{B^{\vee }}$ qui sont la composante en $(P, \o )$ d'un \'el\'ement de $\H (\Theta )$ est un sous-$G\times G$-module de $V_B\otimes V^{\vee }_{B^{\vee }}$.
Le sous-ensemble form\'e des \'el\'ements $\varphi _{P, \so }$ qui v\'erifient par ailleurs la propri\'et\'e $(\P )$ en est un $G\times G$-sous-module. Par ce qui pr\'ec\`ede, ils ont la m\^eme intersection avec $S_P^H(V_B)(\o )\otimes S_{\overline {P}}^H(V^{\vee }_{B^{\vee }})(\o ^{\vee })$ qui est non nulle par 4.2. Il r\'esulte alors de 3.4.2 appliqu\'e au sous-groupe parabolique $P\times \overline {P}$ de $G\times G$ et \`a la repr\'esentation de $M\times M$ d'espace $E_{\so ,B}\otimes E_{\so ^{\vee }, B^{\vee }}\simeq (E_{\so }\otimes E_{\so ^{\vee }})_B$ que ces deux ensembles sont en fait \'egaux, d'o\`u le lemme. \hfill {\fin 2}

\null
\null
{\bf 5.} Pour terminer la preuve de la proposition 0.2, il reste \`a montrer que $\varphi _{P,\so }^{H(\so )}$ a la propri\'et\'e $(\P )$ pour tout $(P, \o)\in \Theta $. Fixons $(P, \o )$. Posons $E=E_{\so }$ et $H=H(\o )$. On notera parfois, pour $\sigma \in \o $, par $E_{\sigma }$ l'espace $E$ s'il est muni de la repr\'esentation $\sigma $.

\null
{\bf 5.1 Lemme:} \it
Soient $\sigma $ dans $\o $, $P'\in \P(M)$ et $v$ un \'el\'ement de $(i_{P'\cap K}^KE)^H$ \`a support dans $(P'\cap K)H$. Soit $a\in A_M$ v\'erifiant les propri\'et\'es de la proposition du num\'ero $2$ relatives \`a  $\overline {P'}$ et $i_{P'}^GE_{\sigma }$.

Alors $(i_{P'}^G\sigma )(1_{HaH})v$ est un \'el\'ement non nul de $S_{\ol {P'}}^H(i_{P'}^GE_{\sigma })$ \`a support dans $(P'\cap K)H$ dont la valeur en $1$ est proportionelle \`a $\sigma (a)v(1)$.

\null
Preuve: \rm
Notons $v_{\sigma }$ l'\'el\'ement de $i_{P'}^GE_{\sigma }$ dont la restriction \`a $K$ est \'egal \`a $v$. \'Ecrivons $P'=MU'$. Pour $k\in K$, on trouve
$$((i_{P'}^G\sigma )(1_{HaH})v)(k)={\mes (HaH)\over \mes (H)}\int _H(i_{P'}^G\sigma (ha)v)(k)dh={\mes (HaH)\over \mes (H)}\int _Hv_{\sigma }(kha)dh.$$
Pour que $kha\in \supp(v_{\sigma })$, il faut que $kh\in P'aHa^{-1}=P'a(H\cap \ol{U'})a^{-1}$. Comme $a$ est positif pour $\ol {P'}$, $a(H\cap \ol{U'})a^{-1}\subseteq H\cap \ol{U'}$, d'o\`u $k\in (P'\cap K)H$. La fonction $ (i_{P'}^G\sigma )(1_{HaH})v$ est donc bien \`a support dans $(P'\cap K)H$.

Sa valeur en $1$ est $\displaystyle {{\mes (HaH)\over \mes (H)}\int _Hv_{\sigma }(ha)dh}$. Apr\`es un calcul \'el\'ementaire qui utilise la d\'ecomposition de $H$ par rapport \`a $P'$ et $\ol {U'}$, on trouve la deuxi\`eme assertion du lemme. 

\hfill {\fin 2}

\null
\null
{\bf 5.2} Fixons une base $\{e_i\}_{i\in I}$ de $E^{H\cap M}$. Notons $\{e_i^{\vee }\}_{i \in I}$ la base duale de $({E^{\vee }})^{H\cap M}$, et $v_i$ (resp. $v_i^{\vee }$) l'\'el\'ement de $(i_{\ol {P}\cap K}^KE)^H$ (resp. $(i_{P\cap K}^KE^{\vee })^H$) \`a support dans $(\ol {P}\cap K)H$ (resp. $(P\cap K)H$) v\'erifiant $v_i(1)=e_i$ (resp. $v_i^{\vee }(1)=e_i^{\vee }$).

\null
{\bf 5.2.1 Lemme:} \it
Soit $\sigma $ un point en position g\'en\'erale de $\o $. Alors $v_i\in S_P^H(i_{\ol {P}}^GE_{\sigma })(\o )$ et de plus, si $w\in W(M,M)$,
$$((J_{P\vert \ol {wP}}(w\sigma )\lambda (w))v_i)(1)
=c_1 \cases {0, &si $w\ne 1$;\cr
             e_i, &si $w=1$,\cr }$$
o\`u $c_1$ est une constante non nulle qui ne d\'epend pas du choix de $\sigma $.

\null
Preuve: \rm 
D'apr\`es le lemme 5.1, $v_i\in S_P^H(i_{\ol {P}}^GE_{\sigma })$. Gr\^ace \`a la proposition 1.3, il suffit alors de prouver la derni\`ere assertion. Soit $w\in W(M,M)$. On a
$$((J_{P\vert \ol {wP}}(w\sigma )\lambda (w))v_i)(1)=
\int _{U\cap wUw^{-1}} v_{i,\sigma }(w^{-1}u)du,$$
si cet op\'erateur d'entrelacement est d\'efini par des int\'egrales convergentes. Par la d\'ecom-position de Bruhat, $w^{-1}u\in \ol {P}H$ seulement si $w=1$. Alors $\int _{U\cap wUw^{-1}} v_{i,\sigma }(u)du=\mes (H\cap U)e_i$. Le cas g\'en\'eral s'en d\'eduit par prolongement analytique. \hfill {\fin 2}

\null
\null
{\bf 5.3 Lemme:} \it
Soit $\sigma \in \o$ en position g\'en\'erale. Les ensembles 
$\{(J_{P\vert \ol{wP}}(\sigma )\lambda (w))$ $v_i\}_{w\in W(M,\so ), i\in I}$ et
$\{(J_{P\vert wP}(\sigma ^{\vee })\lambda (w))v_i^{\vee }\}_{w\in W(M,\so ), i\in I}$ forment des bases de $S_P^H(i_P^GE_{\sigma })$ $(\o )$ et $S_{\ol {P}}^H(i_P^GE^{\vee }_{\sigma ^{\vee }})(\o ^{\vee })$ respectivement. On a
$$<(J_{P\vert \ol{wP}}(\sigma )\lambda (w))v_i, (J_{P\vert w'P}(\sigma ^{\vee })\lambda (w'))v_j^{\vee }>=c_2\delta _{w,w'}\delta _{i,j},$$
o\`u $c_2$ est une constante non nulle qui ne d\'epend pas du choix de $\sigma $.

\null
Preuve: \rm
On va utiliser l'isomorphisme de la proposition 1.3. Soit $w'\in W(M,\o )$. On trouve
$$\eqalign {
(J_{P\vert w'P}(w'\sigma )\lambda (w')J_{P\vert \ol{wP}}(\sigma )\lambda (w)v_i)(1)=&(J_{P\vert w'P}(w'\sigma )J_{w'P\vert \ol{w'wP}}(w'\sigma )\lambda (w'w)v_i)(1)\cr
=&\ j_{P\vert w'P\vert \ol{w'wP}}(w'\sigma )(J_{P\vert \ol{w'wP}}(w'\sigma )\lambda (w'w)v_i)(1).\cr }$$
Gr\^ace \`a 5.2.1, ceci n'est non nul que si $w'=w^{-1}$, et alors $(J_{P\vert \ol{w'wP}}(w'\sigma )\lambda (w'w)v_i)(1)=c_ie_i\in w^{-1}E$. Par ailleurs, 
$j_{P\vert w^{-1}P\vert \ol{P}}(w^{-1}\sigma )=1$. On a donc bien une base de $S_P^H(i_P^GE_{\sigma })$ $(\o )$, puisque son image est une base de $\bigoplus _{w\in W(M,\so )}wE_{\sigma }$.

Soient $w, w'\in W(M,\o )$, $i,j\in I$. Alors
$$\eqalign {
<(J_{P\vert \ol{wP}}(\sigma )\lambda (w))v_i, &(J_{P\vert w'P}(\sigma ^{\vee })\lambda (w'))v_j^{\vee }>=<\lambda (w'^{-1})J_{w'P\vert P}(\sigma )J_{P\vert \ol {wP}}(\sigma )\lambda (w)v_i, v_j^{\vee }>\cr
=&<J_{P\vert w'^{-1}P}(w'^{-1}\sigma )J_{w'^{-1}P\vert \ol {w'^{-1}wP}}(w'^{-1}\sigma )\lambda (w'^{-1}w)v_i, v_j^{\vee }>\cr
=&\ j_{P\vert w'^{-1} \vert \ol{w'^{-1}wP}}(w'^{-1}\sigma )<J_{P\vert \ol {w'^{-1}wP}}(w'^{-1}\sigma )\lambda (w'^{-1}w)v_i, v_j^{\vee }>.\cr }$$
Le support de $v_j^{\vee }$ \'etant $(P\cap K)H$ et l'\'el\'ement \`a gauche \'etant invariant par $H$, ceci n'est non nul que si $(J_{P\vert \ol {w'^{-1}wP}}(w'^{-1}\sigma )\lambda (w'^{-1}w)v_i)(1)\ne 0$. Par 5.2.1, ceci n'est possible que si $w=w'$. Dans ce cas $j_{P\vert w^{-1}P \vert \ol {P}}(w^{-1}\sigma )=1$ et l'expression ci-dessus devient par 5.2.1
$$\int _{M\cap K}\int _{U\cap K}\int _{H\cap \ol {U}} c_1 <(w^{-1}\sigma )(m)e_i, (w^{-1}\sigma ^{\vee })(m) e_j^{\vee }> d\ol{u} du dm
=\mes(H\cap \ol{U}) c_1 \delta _{i,j},$$
d'o\`u le lemme. \hfill {\fin 2}

\null
\null
{\bf 5.4} Posons $\xi _w(\sigma )=c_2^{-1}\sum _{i \in I} v_i\otimes v_i^{\vee }$ pour tout $\sigma \in \o $. On v\'erifie que c'est une application polynomiale sur $\o $. Par ce qui pr\'ec\`ede, on a 
$$\eqalign {
&\sum _{w\in W(M, \so)} (J_{P\vert \ol{wP}}(\sigma )\lambda (w)\otimes J_{P\vert wP}(\sigma ^{\vee })\lambda (w)) \xi _w(w^{-1}\sigma )\cr
=&\sum _{w\in W(M, \so)} c_2^{-1} \sum _{i\in I}(J_{P\vert \ol{wP}}(\sigma )\lambda (w)v_i\otimes J_{P\vert wP}(\sigma ^{\vee })\lambda (w)v_i^{\vee }).\cr }$$
On d\'eduit de 5.3 que ceci est \'egal \`a $\varphi _{P, \so }^H(\sigma )$ pour $\sigma $ en position g\'en\'erale. Par prolongement analytique, ces deux applications sont donc \'egales.

\null
\null
{\bf Lemme:} \it Posons $\xi (\sigma )=\vert W(M,\o )\vert ^{-1}\sum _{w\in W(M,\so )} \xi _w(\sigma )$. Pour tout $\sigma \in \o $, on a
$$\varphi _{P, \so }^H(\sigma )=\sum _{w\in W(M,\so )} (J_{P\vert \overline {wP}}(\sigma )\lambda (w)\otimes J_{P\vert wP}(\sigma ^{\vee })\lambda (w))\ \xi (w^{-1}\sigma ).$$

\null
Preuve: \rm
En utilisant les r\`egles de composition pour les op\'erateurs d'entrelacement, ainsi que les propri\'et\'es de $\varphi _{P, \so }^H$ relatives aux op\'erateurs d'entrelacement, on trouve en effet:

$$\eqalign {
&\sum _{w,w'\in W(M,\so )}(J_{P\vert \ol {wP}}(\sigma )\lambda (w)\otimes J_{P\vert wP}(\sigma ^{\vee })\lambda (w))\xi _{w'}(w^{-1}\sigma ) \cr
=&\sum _{w,w'\in W(M,\so )}(J_{P\vert \ol {ww'P}}(\sigma )\lambda (ww')\otimes J_{P\vert ww'P}(\sigma ^{\vee })\lambda (ww'))\xi _{w'}(w'^{-1}w^{-1}\sigma )\cr
=&\sum _{w\in W(M,\so )}(J_{P\vert wP}(\sigma )\lambda (w)\otimes J_{P\vert wP}(\sigma ^{\vee })\lambda (w))\sum _{w'\in W(M,\so )} j_{P\vert wP\vert \ol {ww'P}}(\sigma )^{-1}\cr
&\ \ \ \ \ j_{P\vert wP\vert ww'P}(\sigma )^{-1}(J_{P\vert \ol {w'P}}(w^{-1}\sigma )\lambda (w')\otimes J_{P\vert w'P}(w^{-1}\sigma ^{\vee })\lambda (w'))\xi _{w'}(w'^{-1}w^{-1}\sigma )\cr
=&\sum _{w\in W(M,\so )} j_{P\vert wP}(\sigma )^{-1}(J_{P\vert wP}(\sigma )\lambda (w)\otimes J_{P\vert wP}(\sigma ^{\vee })\lambda (w))\varphi _{P,\so }^H(w^{-1}\sigma )\cr
=&\ \vert W(M,\o ) \vert \ \varphi _{P,\so }^H(\sigma ).\cr }$$
\hfill {\fin 2}

\null
\null
\null
{\bf R\'ef\'erences bibliographiques:}

\null
[BD] \it Le "centre" de Bernstein \rm - J.N. BERNSTEIN (r\'edig\'e par P.DELIGNE) dans
Repr\'esentations des groupes r\'eductifs sur un corps local, par J.N. Bernstein, P. Deligne,
D. Kazhdan, M.-F. Vign\'eras; Travaux en cours, Hermann, Paris 1984.

\null
[BZ] \it Representations of the group GL$(n,F)$, where $F$ is a local non archimedean field, \rm J. BERNSTEIN et A. ZELEVINSKY, Uspekhi Mat. Nauk. {\bf 31}, 3 (1976), 5-70.

\null
[C] \it Introduction to the theory of admissible representations of $p$-adic reductive groups
\rm - W. CASSELMAN, non publi\'e.

\null
[W] \it La formule de Plancherel pour les groupes $p$-adiques - d'apr\`es Harish-Chandra \rm, J.-L. WALDSPURGER, \`a para\^\i tre.

\bye